\documentclass{article}
\usepackage{amsmath,amssymb,accanthis, mathrsfs,newpxmath, color} 

\usepackage[textsize=tiny]{todonotes}

\newtheorem{thm}{Theorem}[section]

\newtheorem{lem}[thm]{Lemma}
\newtheorem{cor}[thm]{Corollary}
\newtheorem{ex}{Example}[section]
\newtheorem{rem}{Remark}[section]

\title{Backward stochastic difference equations on lattices with application to market equilibrium analysis}
\author{M. Fukasawa$^\sharp$, T. Sato$^{\sharp \flat}$ and J. Sekine$^\sharp$\thanks{Jun Sekine acknowledges Grant-in-Aid for Scientiﬁc Research (C), 23K01450.}\\
{\small  $\sharp$ Graduate School of Engineering Science, Osaka University}\\
{\small $\flat$ Humboldt University of Berlin}
}

\date{}
\begin{document}

\maketitle

\begin{abstract}
    We study backward stochastic difference equations (BS$\Delta$E) 
   driven by a $d$-dimensional stochastic process on a lattice whose increments have only $d+1$ possible values that generates the lattice.
    Regarding the driving process as a $d$ dimensional asset price process,
    we give applications to an optimal investment problem and a market equilibrium analysis, where utility functionals are defined through  BS$\Delta$E. 
\end{abstract}

\section{Introduction}
The theory of Backward Stochastic Differential Equations (BSDE) initiated by Bismut~\cite{Bismut}, Pardoux and Peng~\cite{PP} has been extensively investigated for the last three decades, especially in connection 
with stochastic control, finance and insurance (see e.g.,
\cite{Delong, Zhang}).
Their important applications include dynamic risk measures~\cite{BE} and 
$g$-expectations~\cite{Peng, CHMP}, which extend the notion of expectation and martingale to nonlinear counterparts.
Recent applications to financial economics include \cite{BLE, BFHM, GMS, HMP, MNT, KXZ}. 

While being a powerful theoretical tool, a solution to BSDE is implicit in general and requires discretization in numerical evaluations.
As a discrete analog of BSDE,
Backward Stochastic Difference Equations (BS$\Delta$E)
also have been studied in many papers,
 which fall into two categories.
 The first focuses on the aspect of BS$\Delta$E as a weak approximation of BSDE~\cite{BDM1,BDM2,Nakayama, MPMT, CS,Stadje, BGGL,Tanaka}.
 The second category studies the structure of BS$\Delta$E themselves.
 A general treatment is given in \cite{CE}, while some specific BS$\Delta$E are treated in \cite{CE10,ESC}.
 This paper is in the second category, and deal with a class of BS$\Delta$E including
 the one studied in \cite{Nakayama, Tanaka}, where
 a $d$-dimensional  scaled random walk whose increments have only $d+1$ possible values is employed to approximate a Brownian motion driving a BSDE.
Such a random walk is minimal in a sense among discrete-time processes which have $d$-dimensional Brownian motions as their scaling limits.
This particular BS$\Delta$E is computationally efficient in that it requires to solve only a $d + 1$ dimensional problem while 
the use of the $d$-dimensional symmetric random walk
as in \cite{BDM2} does  a $2^d$ dimensional one.
While a general theory of BS$\Delta$E is given in \cite{CE}, 
in this paper, we are interested in 
the particular structure of such a BS$\Delta$E with the minimality.
 The one-dimensional case 
reduces to BS$\Delta$E on a binomial tree studied in \cite{ESC} in details in terms of dynamic risk measures.
In this paper, we deal with  a general multi-dimensional case.

The $g$-expectation is a part of the solution of a BSDE or BS$\Delta$E that generalizes the expectation and the certainty equivalent of an expected utility. A subclass of them with concavity and a translation invariance property has been employed as the utility functional for market equilibrium analyses in \cite{AD,HPD,CHKP,KXZ}.
In this paper, we also apply our BS$\Delta$E
to a market equilibrium analysis.
In contrast to the preceding studies putting an emphasis on incomplete markets, we are interested in explicit computations in a dynamically complete market.
Anderson and Raimond~\cite{AR} proved the existence of equilibrium in a continuous-time dynamically complete market by means of nonstandard analysis, where an approximation to a Brownian motion by a minimal random walk played a key role.
We consider an even simpler dynamically complete market to derive simple conditions for market equilibrium.
Under a unique equivalent martingale measure, our asset price model is a multi-dimensional extension of the recombining binomial tree.
In our approach, an asset price process is given as a stochastic process taking values 
in a lattice. 
We do not argue the existence of an equilibrium price but characterize the agents' utilities under which the given discrete (in both time and space) price process is in general equilibrium.
This feature is in contrast to the preceding studies \cite{AD,HPD,CHKP} and similar to \cite{HL,S, BFHM,KXZ} in continuous-time.
Our framework includes heterogeneous agents with exponential utilities under heterogeneous beliefs.
Their risk aversion coefficients are allowed to be stochastic and time-varying.
We observe in particular that under equilibrium with heterogeneous beliefs agents trade with each other even when they do not have random endowments to hedge, complementing earlier studies of heterogeneous beliefs~\cite{Varian,HR, Basak,MNT} among others.

Here we give an overview of Section~2.
In Section~2.1, we describe a lattice in $\mathbb{R}^d$ where a stochastic process $\{X_n\}$ takes values, and give some elementary linear algebraic lemmas as a preliminary. 
In Section~2.2, we introduce the process $\{X_n\}$ that is the source of randomness in this paper and generates a filtration.
It is  minimal in the sense that the increment
$\Delta X_n$ is supported on a set $\{v_0,\dots,v_d\}$ of $d+1$ points in $\mathbb{R}^d$.
Some elementary measure change formulas are given also as a preliminary.
In Section~2.3, our BS$\Delta$E
\begin{equation*}
    \Delta Y_n = -g_n(Z_n) + Z_n^\top \Delta X_n,\ \ Y_N = Y
\end{equation*}
is formulated. Since $\{X_n\}$ is minimal,
there exists a unique solution $\{(Y_n,Z_n)\}$ to the above equation without orthogonal martingale terms as in \cite{BDM2, CS}.
The process $\{X_n\}$ itself takes more than $d+1$ points and so, this BS$\Delta$E is different from the one studied in \cite{CE10}.
The $g$-expectation $\mathcal{E}^g_n$ is defined by
$\mathcal{E}^g_n(Y) = Y_n$.
In Section~2.4, we deal with the case $g_n(z) = f_n(X_{n-1},z)$ and $Y_N = h(X_N)$ for deterministic functions $f_n$ and $h$
to derive a nonlinear Feynman-Kac formula, which is a computationally efficient recurrence equation on the lattice for such a deterministic function $u_n$ that
$Y_n = u_n(X_n)$.
In Section~2.5, we state a simple property of the solution regarding a shift in $g_n$, which is not specific to our minimal BS$\Delta$E but given for the readers' convenience.
In Section~2.6, 
we show that the $g$-expectation  is a conditional expectation  when $g_n$ are linear with slope coefficients
included in the convex hull $\Theta$ of the set $\{v_0,\dots,v_d\}$.
This constraint on the slope is a special feature of our BS$\Delta$E.
In Section~2.7, we prove a comparison theorem. A balance condition introduced by \cite{CE} for the comparison theorem to hold is translated in terms of $\Theta$ for our BS$\Delta$E.
In Section~2.8, we prove a robust representation when $g_n$ are concave.
Here again, the set $\Theta$ plays an important role.
In Section~2.9, we show that a translation-invariant filtration-consistent nonlinear expectation is a $g$-expectation.

In Section~3, we regard $\{X_n\}$ as a $d$-dimensional asset price process and consider an optimal investment strategy which maximizes the $g$-expectation of terminal wealth.
By the minimality, the market is complete, extending the well-known binomial tree for a one-dimensional asset.
Our asset price model can be seen as a discrete approximation of the multi-dimensional Bachelier model with constant covariance and general stochastic drift.
An advantage of our use of the minimal process as an approximation is that the completeness of the Bachelier model is preserved. 
Besides, the minimal property naturally arises in a variance swap pricing model as illustrated in Example~\ref{VarSwap}.

In Section~4, we give a market equilibrium analysis.
We consider agents whose utility functionals are $g$ expectations.
Then we seek conditions on those $g$ expectations under which $\{X_n\}$ 
is an equilibrium price process.
In Section~4.1, we formulate the market equilibrium. In Section~4.2, we consider a market with a single agent.
In Section~4.3, we show the existence of a representative agent.
In Section~4.4, we give some explicit computations when agents have similar utility functions with heterogeneous magnitude of risk aversion.
In Section~4.5, we consider agents with heterogeneous beliefs.

Throughout our financial application, we have short maturity problems in our mind and so, for brevity, assume interest rates, dividend rates, and consumption rates to be zero as in \cite{AD,HPD,CHKP, KXZ}.

We use the convention that 
\begin{equation*}
    \sum_{i=m}^n a_i= 0
\end{equation*}
for any sequence $\{a_i\}$ if $m > n$.

\section{BS$\Delta$E on a lattice}
\subsection{Lattice}
We start with describing a lattice.
Let $\{v_1,\dots, v_d\}$ be a basis of $\mathbb{R}^d$. The subset
\begin{equation*}
    L = \left\{\sum_{i=1}^d z_i v_i, \ ; \ z_i \in \mathbb{Z}, i=1,\dots, d\right\}
\end{equation*}
of $\mathbb{R}^d$ is a $d$-dimensional lattice
generated by the basis.
Notice that
\begin{equation}\label{lattice}
    L = \left\{\sum_{i=0}^d n_i v_i, \ ; \ n_i \in \mathbb{N}, i=0,1,\dots, d\right\}
\end{equation}
for $v_0 = -v_1 - \dots - v_d$,
where $\mathbb{N}$ is the set of the nonnegative integers.
Let
\begin{equation*}
    \mathbf{v} = [v_0,v_1,\dots,v_d]
\end{equation*}
be the $d \times (d+1)$ matrix with $v_i$, $i=0,\dots, d$ as its column vectors. Put
\begin{equation*}
    \mathbf{1} = (1,\dots,1)^\top \in \mathbb{R}^{d+1}.
\end{equation*}
The following lemmas will be of repeated use in this paper.
\begin{lem}\label{lem1}
The $(d+1)\times (d+1)$ matrix $(\mathbf{1},\mathbf{v}^\top)$
is invertible.
\end{lem}
{\it Proof: }
We show that the row vectors of $(\mathbf{1},\mathbf{v}^\top)$ are linearly independent.
Suppose
\begin{equation*}
    (\alpha_0,\dots,\alpha_d)(\mathbf{1},\mathbf{v}^\top) = 0
\end{equation*}
or equivalently,
\begin{equation*}
    \sum_{j=0}^d \alpha_j = 0,\ \ \sum_{j=0}^d \alpha_j v_j = 0
\end{equation*}
for scalars $\alpha_j$. By the second equation and the definition of $v_0$, we have
\begin{equation*}
    \sum_{j=1}^d (\alpha_j -\alpha_0) v_j = 0,
\end{equation*}
from which we can conclude $\alpha_j = \alpha_0 $ for all $j$ because $\{v_1,\dots,v_d\}$ is a basis.
Together with the first equation, we then conclude $\alpha_j  = 0$ for all $j$. \hfill{$\square$}

\begin{lem}
Let $y \in \mathbb{R}^{d+1}$.
We have
\begin{equation}\label{syseq}
    y = a(y)\mathbf{1} + \mathbf{v}^\top z(y)
\end{equation}
if and only if
\begin{equation}\label{syseqsol}
    a(y) = \frac{1}{d+1} \mathbf{1}^\top y, \ \ z(y) = (\mathbf{v}\mathbf{v}^\top)^{-1}\mathbf{v}y.
    \end{equation}    
\end{lem}
{\it Proof: }
By Lemma~\ref{lem1}, there uniquely exist $a(y) \in \mathbb{R}$ and
$z(y) \in \mathbb{R}^d$ such that \eqref{syseq} holds.
Since $\mathbf{v}\mathbf{1} = 0$,   from \eqref{syseq},  we have \eqref{syseqsol}. \hfill{$\square$}

\subsection{Probability space}
Let $(\Omega,\mathscr{F},\mathsf{P})$ be a probability space.
For a stochastic process $\{X_n\}_{n \in \mathbb{N}}$, we put
$\Delta X_n = X_n - X_{n-1}$.
Let $\{X_n\}$ be a $d$-dimensional stochastic process
with $\Delta X_n$ taking values in $\{v_0,v_1,\dots,v_d\}$ for all $n \geq 1$ and $X_0 = 0$.
By \eqref{lattice}, $X_n$ takes values in $L$ for all $n \in \mathbb{N}$.
Let $\mathscr{F}_n = \sigma(X_0,\dots,X_n)$, $ n\in \mathbb{N}$,
be the natural filtration generated by $\{X_n\}$, and put
$P_n = (P_{n,0},\dots, P_{n,d})^\top$ for $n\geq 1$, where
\begin{equation*}
    P_{n,j} = \mathsf{P}(\Delta X_n = v_j | \mathscr{F}_{n-1}), \ \ j=0,\dots, d.
\end{equation*}
Note that $\{P_n\}$ is a $\Delta_d$-valued predictable process,
where
\begin{equation*}
    \Delta_d = \{p \in \mathbb{R}^{d+1}\ ; \ \mathbf{1}^\top p = 1, \ \ p \geq 0\}.
\end{equation*}
We assume $P_{n,j}$ is positive for all $n \geq 1$ and $j = 0,\dots, d$. Let $N$ be a natural number ($N \in \mathbb{N}$ and $N \geq 1$).
The following lemma will be of repeated use in this paper.
\begin{lem}\label{lem:mc}
For any $\Delta_d$-valued predictable process $\{\hat{P}_n\}$,
there exists a probability measure $\hat{\mathsf{P}}$ on $(\Omega,\mathscr{F}_N)$ such that
\begin{equation}\label{phatj}
    \hat{P}_{n,j} = \hat{\mathsf{P}}(\Delta X_n = v_j | \mathscr{F}_{n-1}), \ \ j=0,\dots, d, \ \ n=1,\dots, N,
\end{equation}
where $\hat{P}_n = (\hat{P}_{n,0},\dots, \hat{P}_{n,d})^\top$.
\end{lem}
{\it Proof: }
Define $\hat{\mathsf{P}}$ by
    \begin{equation}\label{phat}
        \hat{\mathsf{P}}(A) = \mathsf{E}[L_N1_A],\ \ 
        L_n = \prod_{k=1}^n \left(\sum_{j=0}^d\frac{\hat{P}_{k,j}}{P_{k,j}}1_{\{\Delta X_k = v_j \}}
        \right)
    \end{equation}
for $A \in \mathscr{F}_N$.
The measure $\hat{\mathsf{P}}$ is a probability measure because
\begin{equation*}
    \mathsf{E}\left[   \sum_{j=0}^d\frac{\hat{P}_{n,j}}{P_{n,j}}1_{\{\Delta X_n = v_j \}} \Bigg| \mathscr{F}_{n-1}\right] =   \sum_{j=0}^d\frac{\hat{P}_{n,j}}{P_{n,j}} P_{n,j} = 1
\end{equation*}
and so $L_n$ is a martingale with $L_0 = 1$.
We also have
\begin{equation*}
    \hat{\mathsf{P}}(\Delta X_n =v_j  | \mathscr{F}_{n-1}) = 
    \frac{\mathsf{E}[L_N1_{\{\Delta X_n = v_j\}}|\mathscr{F}_{n-1}]}{\mathsf{E}[L_N |\mathscr{F}_{n-1}]}
    = \mathsf{E}\left[ \frac{\hat{P}_{n,j}}{P_{n,j}}1_{\{\Delta X_n = v_j \}} \Bigg| \mathscr{F}_{n-1}\right] = \hat{P}_{n,j}
\end{equation*}
by the martingale property of $L_n$. \hfill{$\square$}\\

Define a measure $\mathsf{Q}$ on $\mathscr{F}_N$ by
 \begin{equation*}
        \mathsf{Q}(A) = \mathsf{E}[L_N1_A],\ \ 
        L_N = \prod_{n=1}^N \left( \frac{1}{d+1} \sum_{j=0}^d\frac{1}{P_{n,j}}1_{\{\Delta X_n = v_j \}}
        \right).
    \end{equation*}
    Let $\mathsf{E}_\mathsf{Q}$ denote the integration under $\mathsf{Q}$.
    \begin{lem}
    The measure $\mathsf{Q}$ is the unique probability measure on $\mathscr{F}_N$ under which $\{X_n\}$ is a martingale.
        Under $\mathsf{Q}$, 
$\{\Delta X_n\}$ is iid with
\begin{equation*}
    \mathsf{Q}(\Delta X_n = v_j ) = \mathsf{Q}(\Delta X_n =v_j  | \mathscr{F}_{n-1}) = \frac{1}{d+1}
\end{equation*}
for all $n =1,\dots, N$ and $j=0,\dots, d$.
We also have
\begin{equation}\label{meanandvar}
\mathsf{E}_\mathsf{Q}[\Delta X_n |\mathscr{F}_{n-1}]  = 0,\ \ 
 \mathsf{E}_\mathsf{Q}[\Delta X_n (\Delta X_n)^\top |\mathscr{F}_{n-1}] = 
  \frac{1}{d+1} \mathbf{v}\mathbf{v}^\top.
\end{equation}
    \end{lem}
{\it Proof: }
By Lemma~\ref{lem:mc}, $\mathsf{Q}$ is a probability measure with
$ \mathsf{Q}(\Delta X_n =v_j  | \mathscr{F}_{n-1}) = 1/(d+1)$, which implies 
\begin{equation*}
    \mathsf{E}_\mathsf{Q}[\Delta X_n |\mathscr{F}_{n-1}] = \frac{1}{d+1} \mathbf{v}\mathbf{1} = 0.
\end{equation*}
Therefore, $\{X_n\}$ is a martingale with \eqref{meanandvar}.
There is no other such a measure because
\begin{equation*}
    \sum_{j=0}^d \alpha_j = 1, \ \ \sum_{j=0}^d \alpha_j v_j = 0
\end{equation*}
implies $\alpha_j = 1/(d+1)$ as in the proof of Lemma~\ref{lem1}.
Since the conditional law of $\Delta X_n$ given $\mathscr{F}_{n-1}$ is deterministic for every $n$, $\{X_n\}$ is iid. \hfill{$\square$}

\begin{rem}
    \upshape
    For any positive definite $d\times d$ matrix $\Sigma$,
    we can construct such a lattice $L$ that $\mathbf{v}\mathbf{v}^\top = \Sigma$. Indeed, starting with an arbitrary basis, say,
     $\bar{v}_j = e_j$ (the standard basis of $\mathbb{R}^d)$
     with $\bar{v}_0 = -\bar{v}_1 - \dots - \bar{v}_d$ and
     $\bar{\mathbf{v}} = [\bar{v}_0, \dots, \bar{v}_d]$,
     using the Cholesky decomposition  $\Sigma  = CC^\top$ and 
     $\bar{\mathbf{v}}\bar{\mathbf{v}}^\top
     = \bar{C}\bar{C}^\top$,
     define $v_j = C\bar{C}^{-1}\bar{v}_j$, $j=0,\dots, d$.
     Then, $\mathbf{v} = C\bar{C}^{-1}\bar{\mathbf{v}}$ and
     $\mathbf{vv}^\top = C\bar{C}^{-1}\bar{\mathbf{v}}\bar{\mathbf{v}}^\top (\bar{C}^\top)^{-1}C^\top = \Sigma$.
     In particular, we can construct such $v_j$  that $\mathbf{vv}^\top$ is the identity matrix. 
     In this case, a scaling limit of $\{X_n\}$ under $\mathsf{Q}$ is the $d$-dimensional standard Brownian motion.
     Such a set of vectors plays an essential role in proving the existence of continuous-time market equilibrium in Anderson and Raimondo~\cite{AR} by means of nonstandard analysis, where the existence of the vectors are proved in a recursive manner.
     It is also the building block of a $d$-dimensional diamond in
     Topological Crystallography~\cite{Sunada}.
\end{rem}

\subsection{BS$\Delta$E: existence and uniqueness}
 Let $\mathcal{A}$ denote the set of the sequences $g = \{g_n\}_{n=1}^N$
of $\mathscr{F}_{n-1}\otimes \mathscr{B}(\mathbb{R}^d)$ measurable functions $g_n: \Omega \times \mathbb{R}^d \to \mathbb{R}$.
\begin{thm}\label{thm:exist}
Let
$Y_N$ be an $\mathscr{F}_N$-measurable random variable,
 and let $g = \{g_n\} \in \mathcal{A}$.
Then, there exist uniquely an adapted process
$\{Y_n\}_{n=0,\dots,N-1}$ and
an $\mathbb{R}^d$-valued  predictable process $\{Z_n\}_{n=1,\dots,N}$ such that
\begin{equation}\label{bsde}
    \Delta Y_n = -g_n(Z_n) + Z_n^\top \Delta X_n,\ \ n=1,\dots, N.
\end{equation}
\end{thm}
{\it Proof: }
Since $Y_N$ is 
$\mathscr{F}_N$-measurable,
there exists a function $f:L^N \to \mathbb{R}$ such that $Y_N = f(X_1,\dots,X_N)$.
Applying \eqref{syseq} for $y= (y_0,\dots,y_d)^\top$, 
$$y_j = f(X_1,\dots,X_{N-1},X_{N-1} + v_j),$$ we obtain
\eqref{bsde} for $n=N$ with
\begin{equation*}
    Z_N = z(y), \ \ Y_{N-1} = a(y)+g_N(Z_N).
\end{equation*}
It is clear that both $Z_N$ and $Y_{N-1}$ are $\mathscr{F}_{N-1}$-measurable, and such $Z_N$ and $Y_{N-1}$ are unique because $\mathsf{P}_{N,j}$ are positive by the assumption.
By backward induction, we obtain $\{Y_n\}$ and $\{Z_n\}$.
\hfill{$\square$}
\\

For  $g = \{g_n\} \in \mathcal{A}$ fixed, the $\mathscr{F}_n$-measurable random variable $Y_n$ given by Theorem~\ref{thm:exist} is uniquely determined by the $\mathscr{F}_N$-measurable random variable $Y_N$. We write this mapping as
$Y_n = \mathcal{E}^g_n(Y_N)$.
The stochastic process $\{(Y_n,Z_n)\}$ given by Theorem~\ref{thm:exist}
is called the solution of the BS$\Delta$E \eqref{bsde}.

\begin{rem}
    \upshape In the literature, say, in \cite{CE}, BS$\Delta$E is formulated by decomposing $\Delta Y_n$ into a predictable part and a martingale difference part. In our formulation \eqref{bsde}, $\Delta X_n$ is not necessarily a martingale difference. 
    It is a minor reparametrization because 
    \eqref{bsde} can be rewritten as
\begin{equation*}
    \Delta Y_n = - \hat{g}_n(Z_n) + Z_n^\top(\Delta X_n - A_n)
\end{equation*}
with $\hat{g}_n(z) = g_n(z) - z^\top A_n$, $A_n = \mathsf{E}[\Delta X_n | \mathscr{F}_{n-1}]$.
\end{rem}

\begin{ex} \upshape
Let $\gamma > 0$,
$\{(\hat{P}_{n,0},\dots,\hat{P}_{n,d})^\top\}$ be a $\Delta_d$-valued predictable process, 
and
\begin{equation}\label{ex1g}
    g_n(z) = -\frac{1}{\gamma} \log \left(\sum_{j=0}^{d} e^{-\gamma z^\top v_j} \hat{P}_{n,j}\right).
\end{equation}
Then,
\begin{equation}\label{ex1}
    \mathcal{E}^g_n(Y) = -\frac{1}{\gamma} \log 
    \hat{\mathsf{E}}[e^{-\gamma Y}|\mathscr{F}_n], \ \ n=0,1,\dots, N
\end{equation}
for any $\mathscr{F}_N$-measurable random variable $Y$, where
$\hat{\mathsf{E}}$ is the expectation under the measure $\hat{\mathsf{P}}$ on $\mathscr{F}_N$ defined by \eqref{phat}.
    Indeed, by \eqref{syseq},
there exist $\mathscr{F}_{N-1}$-measurable
$A_N$ and $Z_N$ such that
\begin{equation*}
    Y = A_N + Z_N^\top \Delta X_N
\end{equation*}
and so, using Lemma~\ref{lem:mc},
\begin{equation*}
    -\frac{1}{\gamma} \log 
    \hat{\mathsf{E}}[e^{-\gamma Y}|\mathscr{F}_{N-1}]
= A_N +  g_N(Z_N) = Y + g_N(Z_N) - Z_N^\top \Delta X_N = \mathcal{E}^g_{N-1}(Y),
\end{equation*}
which implies \eqref{ex1} for $n=N-1$.
The general case follows by backward induction.
\end{ex}

\begin{thm}\label{thm:mart}
Let $g = \{g_n\} \in \mathcal{A}$ and 
$\{(Y_n,Z_n)\}$ be the solution of the BS$\Delta$E \eqref{bsde}.
Then,
\begin{equation*}
\begin{split}
&     Y_{n-1} = \mathsf{E}_\mathsf{Q}[Y_n|\mathscr{F}_{n-1}] + g_n(Z_n), \ \ 
\\&     Z_n = (d+1) (\mathbf{v}\mathbf{v}^\top)^{-1}  \mathsf{E}_\mathsf{Q}[Y_n \Delta X_n| \mathscr{F}_{n-1}]
=    (\mathbf{v}\mathbf{v}^\top)^{-1}\mathbf{v}\bar{Y}_n,
\end{split}
\end{equation*}
where $\bar{Y}_n = (\bar{Y}_{n,0}, \dots, \bar{Y}_{n,d})^\top$ and
\begin{equation*}
    \bar{Y}_{n,j} = \mathsf{E}_\mathsf{Q}[Y_n|\mathscr{F}_{n-1}, \Delta X_n = v_j]
    = (d+1) \mathsf{E}_\mathsf{Q}[Y_n1_{\{\Delta X_n = v_j\}}| \mathscr{F}_{n-1}].
\end{equation*}
\end{thm}
{\it Proof: }
It follows from \eqref{syseqsol}(or directly, \eqref{bsde}) and \eqref{meanandvar}. \hfill{$\square$}

\subsection{Nonlinear Feynman-Kac formula}\label{FK}
In this subsection we consider the case $g_n(z) = f_n(X_{n-1},z)$ and $Y_N = h(X_N)$ for deterministic functions $f_n : L \times \mathbb{R}^d \to \mathbb{R}$ and $h : L \to \mathbb{R}$.
Let $u_N(x) = h(x)$ and define $u_n : L \to \mathbb{R}$,
$n=0,1,\dots, N-1$
backward inductively by
\begin{equation*}
    u_{n-1}(x) = u_n(x) + \mathcal{L} u_n(x) + f_n(x, (\mathbf{vv}^\top)^{-1}\mathbf{v}\mathcal{N} u_n(x)),
\end{equation*}
where
\begin{equation*}
\begin{split}
   & \mathcal{L} u_n(x)  = \frac{1}{d+1} \sum_{j=0}^d (u_n(x+v_j)-u_n(x)),  
   \\ &  \mathcal{N} u_n(x)
    = (u_n(x+v_0)-u_n(x), \dots, u_n(x+v_d)-u_n(x))^\top.
    \end{split}
\end{equation*}
\begin{thm}
 The unique solution to \eqref{bsde} with $Y_N = h(X_N)$ is given by
 \begin{equation*}
      Y_n = u_n(X_n), \ \ Z_{n+1} = (\mathbf{vv}^\top)^{-1}\mathbf{v}\mathcal{N} u_{n+1}(X_{n}), \ \ n= 0,1,\dots, N-1.
 \end{equation*}
\end{thm}
{\it Proof: }
By definition, $Y_N = h(X_N) = u_N(X_N)$. Suppose $Y_n = u_n(X_n)$. Then,
by Theorem~\ref{thm:mart},
$Z_n = (\mathbf{v}\mathbf{v}^\top)^{-1}\mathbf{v}\bar{Y}_n$,
where
\begin{equation*}
     \bar{Y}_{n,j} = \mathsf{E}_\mathsf{Q}[Y_n|\mathscr{F}_{n-1}, \Delta X_n = v_j] = u_n(X_{n-1} + v_j).
\end{equation*}
Using that $\mathbf{v1} = 0$, we conclude 
$Z_n = (\mathbf{vv}^\top)^{-1}\mathbf{v}\mathcal{N} u_n(X_{n-1})$. Further, again by 
Theorem~\ref{thm:mart},
\begin{equation*}
    Y_{n-1} = 
     \mathsf{E}_\mathsf{Q}[Y_n|\mathscr{F}_{n-1}] + g_n(Z_n)
     = u_n(X_{n-1}) + \mathcal{L} u_n(X_{n-1}) 
     + g_n(Z_n) = u_{n-1}(X_{n-1}),
\end{equation*}
which concludes the proof. \hfill{$\square$}

\subsection{Translation}

\begin{thm}\label{thm:trans}
Let  $g = \{g_n\} \in \mathcal{A}$,
and  $h = \{h_n\}$, $h_n(z) = g_n(z) + B_n$, where
$\{B_n\}$ is a predictable process. Then,
for any $\mathscr{F}_N$-measurable random variable $Y$,
\begin{equation}\label{trans}
    \mathcal{E}^h_\ell(Y) = \mathcal{E}^g_\ell\left( Y + \sum_{i=\ell+1}^NB_i\right),\ \  \ell = 0,\dots,N. 
\end{equation}
\end{thm}
{\it Proof: }
Let $Y^\ast_\ell$ denote the right hand side of \eqref{trans} for $\ell=0,\dots,N$.
Then,
\begin{equation*}
    \Delta Y^\ast_\ell = \mathcal{E}^g_\ell\left( Y + \sum_{i=\ell+1}^NB_i\right)
    -\mathcal{E}^g_{\ell -1}\left( Y + \sum_{i=\ell+1}^NB_i\right) - B_\ell
    = \Delta Y_\ell - B_\ell,
\end{equation*}
where  $\{(Y_n,Z_n)\}$ be the solution of \eqref{bsde} for 
$Y_N = Y + \sum_{i = \ell + 1}^N B_i$.
By Theorem~\ref{thm:mart},
\begin{equation*}
    Z_\ell = (d+1) (\mathbf{v}\mathbf{v}^\top)^{-1}  \mathsf{E}_\mathsf{Q}[Y_\ell \Delta X_\ell| \mathscr{F}_{\ell-1}]
    = (d+1) (\mathbf{v}\mathbf{v}^\top)^{-1}  \mathsf{E}_\mathsf{Q}[Y^\ast_\ell \Delta X_\ell| \mathscr{F}_{\ell-1}] =: Z^\ast_\ell.
\end{equation*}
 Therefore,
\begin{equation*}
    \Delta Y^\ast_\ell 
    = \Delta Y_\ell - B_\ell
    = -g_\ell(Z_\ell) +  Z_\ell^\top \Delta X_\ell -B_\ell 
    = -h_\ell(Z^\ast_\ell) + (Z^\ast_\ell)^\top \Delta X_\ell
\end{equation*}
for each $\ell = 0,\dots, N$, that is, $\{(Y^\ast_n,Z^\ast_n)\}$ solves
\eqref{bsde} with $g_n$ replaced by $h_n$.
In particular, $Y^\ast_\ell = \mathcal{E}^h_\ell(Y)$.
\hfill{$\square$}

\subsection{Linear BS$\Delta$E}
Let $\Theta$ be the closed convex hull spanned by $\{v_0,v_1,\dots,v_d\}$,
or equivalently,
\begin{equation*}
    \Theta = \{\mathbf{v}p \ ; \ p \in \Delta_d\}.
\end{equation*}
\begin{ex}
\upshape
The triangular lattice of $\mathbb{R}^2$ is generated by
\begin{equation*}
    v_1 =\frac{1}{\sqrt{6}} \begin{pmatrix}
    0\\
-2
    \end{pmatrix}, \ \ 
    v_2 = \frac{1}{\sqrt{6}}
    \begin{pmatrix}
        \sqrt{3} \\
        1
    \end{pmatrix}.
\end{equation*}
In this case, $\mathbf{v}\mathbf{v}^\top = I$
and $\Theta$ is an equilateral triangle.
\end{ex}

\begin{thm}\label{thm:linear}
    Let $g_n(z) = A_n^\top z + B_n$ for 
    a $\Theta$-valued predictable process $\{A_n\}$ and 
  a predictable process $\{B_n\}$,
    $n=1,\dots,N$. Then
    \begin{equation*}
        \mathcal{E}^g_n (Y)  = \hat{\mathsf{E}}\left[Y
   +\sum_{i=n+1}^N B_i     
    \Bigg|\mathscr{F}_n\right], \ \ n=0,1,\dots,N
    \end{equation*}
    for any $\mathcal{F}_N$-measurable random variable $Y$,
    where $\hat{\mathsf{E}}$ is the expectation under the measure $\hat{\mathsf{P}}$ on $\mathscr{F}_N$ defined by \eqref{phat}
    with $\hat{P}_n = (\hat{P}_{n,0},\dots, \hat{P}_{n,d})^\top$ such that $A_n = \mathbf{v}\hat{P}_n$.
\end{thm}
{\it Proof: }
By Lemma~\ref{lem:mc},
\begin{equation*}
    \hat{\mathsf{E}}[\Delta X_n |\mathcal{F}_{n-1}] = 
    \frac{\mathsf{E}[L_N \Delta X_n | \mathcal{F}_{n-1}]}{\mathsf{E}[L_N | \mathcal{F}_{n-1}]}
= \sum_{j=0}^d v_j \hat{P}_{n,j}
= \mathbf{v}\hat{P}_n = A_n
\end{equation*}
for all $n$.
On the other hand, from \eqref{bsde}, we have
\begin{equation*}
    Y= Y_N = Y_n + \sum_{i=n+1}^N (-g_i(Z_i) + Z_i^\top \Delta X_i) =  Y_n + \sum_{i=n+1}^N (-B_i + Z_i^{\top}(\Delta X_i - A_i)).
\end{equation*}
Taking the conditional expectation under $\hat{\mathsf{P}}$, we get the conclusion.
\hfill{$\square$}

\begin{ex}\upshape
Let $N=1$, $d=1$, $\Omega = \{+,-\}$, $v_0 = -1, v_1 = 1$,
and $\Delta X_1(\pm) = \pm 1$.
Then, $L = \mathbb{Z}$, $\mathbf{v} = (-1,1)$, $\Theta = [-1,1]$ and
$\mathbf{vv}^\top = 2$. Following the proof of Theorem~\ref{thm:exist},
the solution of the linear BS$\Delta$E
$\Delta Y_1 = -AZ_1 + Z_1\Delta X_1$ can be constructed as
\begin{equation*}
    Y_0 = \frac{Y_1(+) + Y_1(-)}{2} + A \frac{Y_1(+)-Y_-(-)}{2} = 
    \frac{1+A}{2} Y_1(+) + \frac{1-A}{2} Y_1(-)
\end{equation*}
for any $A \in \mathbb{R}$.
The expression given by Theorem~\ref{thm:linear} is 
$Y_0 = \hat{E}[Y_1]$, which can be directly seen
with $\hat{P}_1 = ((1-A)/2, (1+A)/2)^\top$ for $A \in \Theta = [-1,1]$.
For $A \notin \Theta$, we observe that $Y_0$ is not increasing in 
either $Y_1(+)$ or $Y_1(-)$.
\end{ex}

\subsection{Comparison theorem}
 Let $\mathcal{B}$ denote the set of the sequence 
$g = \{g_n\}_{n=1}^N \in \mathcal{A}$ with
 \begin{equation}\label{balance}
        g_n(z_2) - g_n(z_1) \geq \min_{\theta\in\Theta}\theta^\top(z_2-z_1)
    \end{equation}
    for all $z_1, z_2 \in \mathbb{R}^d$.
\begin{thm}\label{thm:comp}
For $i=1,2$,
let  $Y^{(i)}$  be $\mathscr{F}_N$-measurable random variables with $Y^{(1)} \geq Y^{(2)}$, and  
$g^{(i)}  = \{g^{(i)}_n\} \in \mathcal{A}$ with $g^{(1)}_n \geq g^{(2)}_n$.
Let
$\mathcal{E}^{(i)}_n(Y^{(i)})$ denote $\mathcal{E}^g_n(Y^{(i)})$ 
for $g = g^{(i)}$, $i=1,2$ respectively.
Assume also $g^{(i)} \in \mathcal{B}$ for either $i=1$ or $i=2$.
Then,
\begin{equation*}
    \mathcal{E}^{(1)}_n(Y^{(1)}) \geq  \mathcal{E}^{(2)}_n(Y^{(2)}), \ \
    n=0,1,\dots, N.
\end{equation*}
\end{thm}
{\it Proof: }
Note first that
\begin{equation}\label{bal}
    \min_{\omega \in \Omega} z^\top \Delta X_n(\omega) 
    = \min_{p \in \Delta_d} z^\top \mathbf{v}p = \min_{\theta \in \Theta}z^\top \theta
\end{equation}
for any $z \in \mathbb{R}^d$ and $n$.
Therefore, under \eqref{balance} for $g= g^{(i)}$,  $g^{(i)}$ is balanced in the terminology of
\cite{CE}, and so the result follows from Theorem~3.2 of \cite{CE}.
Here we repeat essentially the same proof for the readers' convenience.
Let $\{(Y^{(i)}_n,Z^{(i)}_n)\}$ be the solution of \eqref{bsde}
with $g = g^{(i)}$ and $Y_N = Y^{(i)}$.
We have $Y^{(1)}_N \geq Y^{(2)}_N$ by assumption.
Suppose $Y^{(1)}_k \geq Y^{(2)}_k$ for some $k$. Then,
\begin{equation*}
     0 \leq          Y^{(1)}_k - Y^{(2)}_k =  Y^{(1)}_{k-1} - Y^{(2)}_{k-1}
          - g^{(1)}_k(Z^{(1)}_k)
         + g^{(2)}_k(Z^{(2)}_k) + (Z^{(1)}_k-
         Z^{(2)}_k)^\top \Delta X_k
\end{equation*}
and so
\begin{equation*}
    (Z^{(1)}_k-
         Z^{(2)}_k)^\top \Delta X_k 
         \geq  - Y^{(1)}_{k-1} + Y^{(2)}_{k-1}
         + g^{(1)}_k(Z^{(1)}_k)
         - g^{(2)}_k(Z^{(2)}_k).
\end{equation*}
Since the right hand side is $\mathscr{F}_{k-1}$-measurable, this implies further
\begin{equation*}
   \min_{\theta \in \Theta} \theta^\top (Z^{(1)}_k-
         Z^{(2)}_k) \geq  - Y^{(1)}_{k-1} + Y^{(2)}_{k-1}
         + g^{(1)}_k(Z^{(1)}_k)
         - g^{(2)}_k(Z^{(2)}_k).
\end{equation*}
by \eqref{bal}. Therefore,
 \begin{equation*}
     \begin{split}
        Y^{(1)}_{k-1} - Y^{(2)}_{k-1} & \geq  
         g^{(1)}_k(Z^{(1)}_k)
         - g^{(2)}_k(Z^{(2)}_k) -  \min_{\theta \in \Theta} \theta^\top (Z^{(1)}_k-
         Z^{(2)}_k) 
        \\ & = g^{(1)}_k(Z^{(2)}_k)
         - g^{(2)}_k(Z^{(2)}_k) +
          g^{(1)}_k(Z^{(1)}_k)
         - g^{(1)}_k(Z^{(2)}_k) - 
            \min_{\theta \in \Theta} \theta^\top (Z^{(1)}_k-
         Z^{(2)}_k) 
         \\ & \geq 0
     \end{split}
 \end{equation*}
under \eqref{balance} for $g_n = g^{(1)}_n$ and also
\begin{equation*}
     \begin{split}
        Y^{(1)}_{k-1} - Y^{(2)}_{k-1} & \geq  
         g^{(1)}_k(Z^{(1)}_k)
         - g^{(2)}_k(Z^{(2)}_k) -  \min_{\theta \in \Theta} \theta^\top (Z^{(1)}_k-
         Z^{(2)}_k) 
        \\ & = g^{(1)}_k(Z^{(1)}_k)
         - g^{(2)}_k(Z^{(1)}_k) +
          g^{(2)}_k(Z^{(1)}_k)
         - g^{(2)}_k(Z^{(2)}_k) - 
            \min_{\theta \in \Theta} \theta^\top (Z^{(1)}_k-
         Z^{(2)}_k) 
         \\ & \geq 0
     \end{split}
 \end{equation*}
 under \eqref{balance} for $g_n = g^{(2)}_n$.
The result then follows by induction.
\hfill{$\square$}

\begin{rem}\label{rem:comp} \upshape
    A sufficient condition for $g_n$ to meet \eqref{balance} is that $g_n(z)$ is continuously differentiable in $z$ with $\nabla g_n(z)$ taking values in $\Theta$.
    Indeed, by Taylor's theorem,
    \begin{equation*}
    g_n(z_1) - g_n(z_2) = A_n^\top (z_1-z_2),\ \ 
    A_n = \int_0^1 \nabla g_n(z_2 + t(z_1 -z_2))\,\mathrm{d}t
\end{equation*}
and then notice that $A_n$ is $\Theta$-valued because
 $\Theta$ is a convex set.
\end{rem}
\begin{ex}[ Locally entropic monetary utility]
    \upshape 
    Let $\{(\hat{P}_{n,0},\dots,\hat{P}_{n,d})^\top\}$ be a $\Delta_d$-valued predictable process, 
    $\{B_n\}$ and $\{G_n\}$ be positive predictable processes,
and
\begin{equation}\label{ex3g}
    g_n(z) = -\frac{1}{G_n} \log \left(\sum_{j=0}^{d} e^{-G_n z^\top v_j} \hat{P}_{n,j}\right) - \frac{1}{G_n}\log B_n.
\end{equation}
Using a similar calculation, shown in Example 2.1, we deduce the relation
\[
Y_{n-1}=-\frac{1}{G_n} 
\left\{
\log \hat{\mathsf{E}}\left[ e^{-G_n Y_n} \middle| {\mathcal F}_{n-1}\right]
+\log B_n
\right\},
\quad n=N,\dots, 1.
\]
In particular, when $B_n=1$, 
 ${\mathcal E}_{n-1}^g$ is 
  locally
the minus of the entropic risk measure with  risk aversion parameter $G_n$
extending  \eqref{ex1}.
In contrast to the dynamic entropic risk measure studied in \cite{AP},
we have the 
time-consistency property
$\mathcal{E}_m^g(\mathcal{E}_n^g(Y))= 
\mathcal{E}_m^g(Y)$ for any $m \leq n$
 when $B_n = 1$ for all $n$
even if the process $\{G_n\}$ is not constant.
We allow $B_n \neq 1$ in order to include an example in Section~4.
We call $\mathcal{E}^g_n$ a  locally entropic monetary utility.
We have
    \begin{equation}\label{ngn}
        \nabla g_n(z) = \frac{\hat{\mathsf{E}}[\Delta X_n e^{-G_n z^\top \Delta X_n}|\mathscr{F}_{n-1}]}{\hat{\mathsf{E}}[e^{-G_n z^\top \Delta X_n}|\mathscr{F}_{n-1}]} = \mathbf{v}\hat{P}_n(z),
    \end{equation}
    where $\hat{P}_n(z) = (\hat{P}_{n,0}(z),\dots,\hat{P}_{n,d}(z))^\top$ and 
    \begin{equation}\label{pnz}
        \hat{P}_{n,j}(z)= \frac{e^{-G_n z^\top v_j} \hat{P}_{n,j}}{\sum_{k=0}^d e^{-G_n z^\top v_k} \hat{P}_{n,k}}.
    \end{equation}
    Since $\hat{P}_n(z)$ is continuous in $z$ and $\Delta_d$-valued for all $n$, by Remark~\ref{rem:comp},
 the assumptions of Theorem~\ref{thm:comp} on $g^{(i)}_n$ are satisfied.
\end{ex}


\subsection{Concavity}
Let $\mathcal{C}$ denote the set of $g = \{g_n\} \in \mathcal{B}$
with $g_n(z)$ being concave in $z$ for all $n$.

\begin{lem}
    Let $g = \{g_n\} \in \mathcal{A}$.
    Then, $g \in \mathcal{C}$ if and only if
    \begin{equation}\label{legendre}
        g_n(z) = \min_{\theta \in \Theta} \left\{z^\top \theta + b_n(\theta)\right\},
    \end{equation}
    where
    \begin{equation}\label{bnt}
        b_n(\theta) = \sup_{z \in \mathbb{R}^d} \left\{g_n(z)-z^\top \theta\right\}.
    \end{equation}
\end{lem}
{\it Proof: }
If $g_n$ is concave, then it is continuous on the interior of its domain that is $\mathbb{R}^d$. Therefore by a well-known fact on the Legendre transform,
we have
 \begin{equation*}
        g_n(z) = \inf_{x \in \mathbb{R}^d} \left\{z^\top x + b_n(x)\right\}.
        \end{equation*}
    Let $x \notin \Theta$.
Since $\Theta$ is a closed convex set of $\mathbb{R}^d$, by the Hahn-Banach theorem, there exists
$z_0 \in \mathbb{R}^d$ such that
\begin{equation*}
    \min_{\theta \in \Theta}z_0^\top \theta  > z_0^\top x.
\end{equation*}
Using \eqref{balance}, for $z = \alpha z_0$, $\alpha > 0$.
\begin{equation*}
    g_n(z)  - z^\top x \geq g_n(0) + \min_{\theta \in \Theta} \theta^\top z -z^\top x
    = g_n(0)+ \alpha \min_{\theta \in \Theta} z_0^\top(\theta -x).
\end{equation*}
Since the last term is positive, letting $\alpha \to \infty$,
we conclude $b_n(x)=\infty$.  This implies
\begin{equation*}
        g_n(z) = \inf_{\theta \in \Theta} \left\{z^\top \theta + b_n(\theta)\right\}
        = \inf_{(\theta,b) \in A_n} \left\{z^\top \theta + b\right\},
    \end{equation*}
    where 
    $A_n = \{(\theta,b) \in \Theta; g_n(w) \leq w^\top \theta + b \text{ for all } w  \in \mathbb{R}^d\}$.
    Fix $n$ and $z$ and then, take a sequence $\{(\theta_k,b_k)\} \subset A_n$ such that
    $z^\top \theta_k + b_k \to g_n(z)$. Since $\Theta$ is compact,
    there exists a converging subsequence $\{\theta_{k_j}\}$ with limit $\theta_\ast \in \Theta$. We have $b_{k_j} = z^\top\theta_{k_j} + b_{k_j} - z^\top \theta_{k_j} \to g_n(z) - z^\top \theta_\ast =:b_\ast$. Also,
    $w^\top \theta_{k_j} + b_{k_j} \geq g_n(w)$ for all $w$ 
    implies $w^\top \theta_\ast + b_\ast \geq g_n(w)$ for all $w$, hence
    $(\theta_\ast,b_\ast) \in A_n$.
    Thus we obtain \eqref{legendre}.
    Conversely, if \eqref{legendre} is true, then
    $g_n(z)$ is concave as being the minimum of concave (affine) functions.
    Since
    \begin{equation*}
   z_1^\top \theta + b_n(\theta)  
   =   z_2^\top \theta + b_n(\theta) + 
   \theta^\top(z_1-z_2) \geq 
  z_2^\top \theta + b_n(\theta) +  \min_{\theta \in \Theta} \theta^\top(z_1-z_2),
    \end{equation*}
    we derive \eqref{balance} from \eqref{legendre}.
    \hfill{$\square$}\\



Let $\mathcal{P}$ denote the set of the probability measures 
on $(\Omega,\mathscr{F}_N)$ absolutely continuous with respect to $\mathsf{P}$.
For $\hat{\mathsf{P}} \in \mathcal{P}$,
there corresponds a $\Delta_d$-valued predictable process
$\{\hat{P}_n\}$ by \eqref{phatj}. 
The measure $\hat{\mathsf{P}}$ is recovered from $\{\hat{P}_n\}$ by \eqref{phat}. 
Let $\hat{\mathsf{E}}$ denote the expectation under $\hat{\mathsf{P}}$ and
\begin{equation*}
    c^g_n(\hat{\mathsf{P}}) = \hat{\mathsf{E}}\left[\sum_{i=n+1}^N b_i(\mathbf{v}\hat{P}_i)\Bigg|\mathscr{F}_n\right],
\end{equation*}
where $b_i$ is associated with $g  = \{g_n\} \in \mathcal{C}$ via \eqref{legendre}.
\begin{thm}\label{thm:robust}
Let $g  \in \mathcal{C}$. Then, for any $\mathscr{F}_N$-measurable random variable $Y$,
\begin{equation}\label{robust}
    \mathcal{E}^g_n(Y) = \min_{\hat{\mathsf{P}}\in \mathcal{P}} \left\{
    \hat{\mathsf{E}}[Y|\mathscr{F}_n] + c^g_n(\hat{P}) \right\}, \ \ n=0,1,\dots, N.
\end{equation}
then it does  in \eqref{robust}.
\end{thm}
{\it Proof: }
By \eqref{legendre}, we have
\begin{equation*}
    g_n(z) \leq z^\top \mathbf{v}\hat{P}_n + b_n(\mathbf{v}\hat{P}_n)
\end{equation*}
for any $\hat{P} \in \mathcal{P}$.
Therefore, by Theorems~\ref{thm:linear} and \ref{thm:comp}, we have
\begin{equation*}
    \mathcal{E}^g_n(Y) \leq  \hat{\mathsf{E}}[Y|\mathscr{F}_n] + c^g_n(\hat{P})
\end{equation*}
for any $\hat{P} \in \mathcal{P}$.
On the other hand, for any $Y \in \mathscr{F}_N$,
there exists the solution $\{(Y_n,Z_n)\}$ of \eqref{bsde} with $Y_N= Y$.
By \eqref{legendre}, 
 there exists $p = p(Z_n)$ such that
\begin{equation*}
    g_n(Z_n)  =  Z_n^\top \mathbf{v}p +  b_n(\mathbf{v}p).
\end{equation*}
Let $\hat{P} \in \mathcal{P}$ be associated with $\hat{P}_n = p(Z_n)$.
Then, by Theorems~\ref{thm:trans} and \ref{thm:linear},
\begin{equation*}
    \mathcal{E}^g_n(Y) =  \hat{\mathsf{E}}[Y|\mathscr{F}_n] + c^g_n(\hat{P}),
\end{equation*}
which implies \eqref{robust}. 
\hfill{$\square$}

\begin{cor}
Let $g \in \mathcal{C}$. Let $Y$ and $Y^\prime$ be
$\mathscr{F}_N$-measurable random variables.
\begin{enumerate}
    \item If $Y \geq Y^\prime$, then $\mathcal{E}^g_n(Y) \geq \mathcal{E}^g_n(Y^\prime)$, $n=0,1,\dots, N$.
    \item For any $\mathscr{F}_n$-measurable $[0,1]$-valued random variable $\lambda$,
     \begin{equation*}
       \mathcal{E}^g_n(\lambda Y + (1-\lambda)Y^\prime)
        \geq \lambda
          \mathcal{E}^g_n(Y) + 
           (1-\lambda) \mathcal{E}^g_n(Y^\prime), \ \ 
           n=0,1,\dots, N.
    \end{equation*}
\end{enumerate}
\end{cor}
%

\begin{ex}
\upshape
 Let $\Theta_n \subset \Theta$  and
\begin{equation*}
    g_n(z) = \inf_{\theta \in \Theta_n} z^\top \theta, \ \ 
    n= 1,\dots, N.
\end{equation*}
Then, we have \eqref{legendre} with $b_n$ such that
 $b_n(\theta) = 0$ for if $\theta \in \bar{\Theta}_n$ while
 $b_n(\theta) = \infty$ otherwise, where $\bar{\Theta}_n$ is the closure of 
 $\Theta_n$.
In particular when $\Theta_n = \Theta$, by Theorem~\ref{thm:robust},
\begin{equation*}
    \mathcal{E}^g_n(Y) = \min_{\hat{\mathsf{P}} \in \mathcal{P}} \hat{\mathsf{E}}[Y |\mathscr{F}_n],
     \ \ 
    n= 1,\dots, N
\end{equation*}
for any $\mathscr{F}_N$-measurable random variable $Y$.
Note also that
\begin{equation*}
    \mathcal{E}^g_0(Y) = \min_{\hat{\mathsf{P}} \in \mathcal{P}} \hat{\mathsf{E}}[Y] = \min_{\omega \in \Omega} Y(\omega).
\end{equation*}
When $\Theta_n = \{\mathbf{v}P_n^{(1)}, \dots, \mathbf{v}P^{(m)}_n\}$ for a $\Delta_d$-valued predictable process $\{P^{(i)}_n\}$,
letting $\mathsf{E}^{(x)}$ denote the expectation under 
the measure determined by $\{P^{(x_n)}_n\}$ for $x = (x_1,\dots,x_N) \in \{1,\dots,m\}^N$,
 by Theorem~\ref{thm:robust}, we have
\begin{equation*}
    \mathcal{E}^g_n(Y) = \min_x \mathsf{E}^{(x)}[Y |\mathscr{F}_n],
     \ \ 
    n= 1,\dots, N.
\end{equation*}
\end{ex}

\subsection{Filtration consistent nonlinear expectations}
Inspired by Coquet et al.~\cite{CHMP},
we call $\mathcal{E}:L^0(\Omega,\mathscr{F}_N,\mathsf{P}) \to \mathbb{R}$
a filtration consistent nonlinear expectation if
\begin{enumerate}
    \item $Y \geq Y^\prime \Rightarrow \mathcal{E}(Y) \geq \mathcal{E}(Y^\prime)$,
    \item $Y \geq Y^\prime$ and $ \mathcal{E}(Y) = \mathcal{E}(Y^\prime) 
    \Rightarrow Y = Y^\prime$, 
    \item $\mathcal{E}(c) = c$ for any constant $c \in \mathbb{R}$, and
    \item for any $n=1,\dots,N$ and  $Y$,
    there exists an $\mathscr{F}_n$-measurable $\eta$ such that 
    $\mathcal{E}(Y1_A) = \mathcal{E}(\eta 1_A)$
    for any $A \in \mathscr{F}_n$.
\end{enumerate}
Such $\eta$ is uniquely determined as shown in \cite{CHMP}.
Let it be denoted by $\mathcal{E}_n(Y)$.
It follows that
\begin{equation}\label{tower}
    \mathcal{E}_n(\mathcal{E}_m(Y)) = \mathcal{E}_n(Y)
\end{equation}
for any $m \geq n$ and 
\begin{equation}\label{local}
    \mathcal{E}_n(Y)1_A = \mathcal{E}_n(Y1_A)
\end{equation}
for any $A \in \mathscr{F}_n$.

\begin{thm}
Let $g = \{g_n \} \in \mathcal{A}$ and
assume that for $n = 1,\dots, N$,
\begin{enumerate}
    \item  $g_n(0) = 0$ and
    \item  
for any $z_1, z_2 \in \mathbb{R}^d$,
    \begin{equation*}
    g_n(z_1) -g_n(z_2) \geq \min_{\theta \in \Theta} \theta^\top (z_1-z_2)
\end{equation*}
with equality holding only if $z_1 = z_2$.
\end{enumerate}
Then, $\mathcal{E}^g_0$ is a filtration consistent nonlinear expectation 
with $\mathcal{E}_n = \mathcal{E}^g_n$
and 
a translation invariance property:
\begin{equation}\label{trin}
    \mathcal{E}_n(Y + \eta) = \mathcal{E}_n(Y) + \eta
\end{equation}
for any $\mathscr{F}_N$-measurable random variable $Y$,
$\mathscr{F}_n$-measurable random variable $\eta$
and   $n=1,\dots, N$.
\end{thm}
{\it Proof: } 
From Theorem~\ref{thm:comp} and its proof, we observe the first two properties of filtration consistent nonlinear expectation.
By $g_n(0) = 0$ we derive \eqref{tower} and \eqref{local}, from which the other properties follow.
\hfill{$\square$}
\\

The following theorem is a discrete analog of Theorem~7.1 of \cite{CHMP}.
\begin{thm}
Let $\mathcal{E}$ be a filtration consistent nonlinear expectation with the translation invariance property \eqref{trin}.
Let $g_n(z) = \mathcal{E}_{n-1}(z^\top \Delta X_n)$.
Then, 
$\mathcal{E}_n = \mathcal{E}^g_n$.
\end{thm}
{\it Proof: } We have $\mathcal{E}_N(Y) = Y$ for any $Y$,
so the claim is true for $n=N$.
Assume $\mathcal{E}_k(Y) = \mathcal{E}^g_k(Y)$ for $k \geq n$. Then,
by \eqref{syseq}, there exists $\mathscr{F}_{n-1}$-measurable $A_n$ and $Z_n$ such that
\begin{equation*}
    \mathcal{E}_n(Y) = A_n  + Z_n^\top \Delta X_n.
\end{equation*}
Then, by \eqref{tower} and \eqref{trin},
\begin{equation*}
      \mathcal{E}_{n-1}(Y) = A_n  + \mathcal{E}_{n-1}(Z_n^\top \Delta X_n).
\end{equation*}
The last term is $g_n(Z_n)$ by \eqref{local}.
Therefore,
\begin{equation*}
     \mathcal{E}_n(Y)-
      \mathcal{E}_{n-1}(Y) = - g_n(Z_n) + Z_n^\top \Delta X_n,
\end{equation*}
which implies that $\mathcal{E}_{n-1}(Y) = \mathcal{E}_{n-1}^g(Y)$.
The result follows by induction.
\hfill{$\square$}

\section{Monetary utility maximization}
Now we consider $\{X_n\}$ to be a $d$-dimensional asset price process.
The lattice $L$ is then understood as the price grid.
For any $\{\mathscr{F}_n\}$ predictable $\mathbb{R}^d$-valued process $Z = \{Z_n\}$
and $w \in \mathbb{R}$,
\begin{equation*}
    W_n(w,Z) := w + \sum_{i=1}^n Z_i^\top \Delta X_i
\end{equation*}
represents the wealth process associated with the portfolio strategy $Z$ and the initial wealth $w$.
Applying Theorem~\ref{thm:exist} with $g_n = 0$, we observe that the market is complete.

Consider an agent whose utility functional is
$\mathcal{E}^g_0$ for $g \in \mathcal{C}$.
This utility is monetary in the sense that
for all $n = 0, \dots, N$,
$\mathcal{E}_n^g(Y+A) = \mathcal{E}_n^g(Y)+A$
for any $\mathscr{F}_N$-measurable random variable $Y$ and 
$\mathscr{F}_n$-measurable random variable $A$.
When assuming $g_n(0) = 0$ for all $n$ in addition, the utility is normalized in the sense that $\mathcal{E}^g_n(0)=0$, $n =0,\dots, N$,
and it is time-consistent in the sense that
$\mathcal{E}_m^g(\mathcal{E}_n^g(Y))= 
\mathcal{E}_m^g(Y)$ for any $m \leq n$
and for any $\mathscr{F}_N$-measurable random variable $Y$.
The simplest example is 
$\mathcal{E}^g_n(Y) = \mathsf{E}[Y|\mathscr{F}_n]$
corresponding to $g_n(z) = z^\top \mathsf{E}[\Delta X_n |\mathscr{F}_{n-1}]$.
More generally, $\mathcal{E}^g_n(Y)$ is a conditional expectation with respect to a probability measure when $g_n$ are linear.
The driver $\{g_n\}$ should reflect
the agent's belief on the distribution of the price process $\{X_n\}$.
For example, if $g_n = 0$, then $\mathcal{E}^g_n(Y) = \mathsf{E}_\mathsf{Q}[Y|\mathscr{F}_n]$ irrespectively of $\mathsf{P}$.
The choice of nonlinear $\{g_n\}$ accommodates a nonlinear evaluation of risk, extending the exponential utility \eqref{ex1}.
In light of Theorem~\ref{thm:robust},
our problem can be interpreted as a robust utility maximization.

The agent's objective is to maximize
$\mathcal{E}^g_0(H + W_N(w,\pi))$ among the predictable process $\pi = \{\pi_n\}$, where $H$ is a given $\mathscr{F}_N$-measurable random variable representing an initial endowment of the agent. 
Since the utility is monetary, it suffices to treat the case $w = 0$.
Let 
\begin{equation*}
    Y^\pi_n =\mathcal{E}^g_n\left(H + W_N(0,\pi)- W_n(0,\pi) \right) =\mathcal{E}^g_n(H + W_N(0,\pi))- W_n(0,\pi).
\end{equation*}
Then the problem is equivalent to maximizing $Y^\pi_0$ among $\pi$.

\begin{thm}\label{thm:max}
Assume 
that there exists 
a predictable process $\{Z^\dagger_n\}$
 such that
\begin{equation}\label{dagger}
    g_n(Z^\dagger_n) = \sup_{z \in \mathbb{R}^d} g_n(z), \ \ n=1,\dots, N.
\end{equation}
Then,
\begin{equation*}
  \max_{\pi} Y^\pi_0 = \max_{\pi}\mathcal{E}^g_0(H + W_{N} (0,\pi))
    = \mathcal{E}^g_0(H + W_{N}(0,\pi^\ast)),
\end{equation*}
where 
$\pi^\ast_n = Z^\dagger_n -  Z^\ast_n$,
$ Z^\ast_n = Z^H_n + Z^g_n$ and
\begin{equation}\label{rep}
\begin{split}
& Z^g_n = \sum_{i=n+1}^N
(d+1) (\mathbf{v}\mathbf{v}^\top)^{-1}\mathsf{E}_\mathsf{Q}[g_i(Z_i^\dagger) \Delta X_n | \mathscr{F}_{n-1}],\\
& Z^H_n = 
(d+1) (\mathbf{v}\mathbf{v}^\top)^{-1}\mathsf{E}_\mathsf{Q}[H\Delta X_n | \mathscr{F}_{n-1}]
\end{split}
\end{equation}
for $n=1,\dots, N$. Moreover,
\begin{equation}\label{yast}
     Y^{\pi^\ast}_n = \mathsf{E}_\mathsf{Q}[
    H  | \mathscr{F}_n]
     + \sum_{i=n+1}^N \mathsf{E}_\mathsf{Q}[  g_i(Z_i^\dagger) |  \mathscr{F}_n]
\end{equation}
for $n=0,1,\dots, N$.
If in addition $Z^\dagger_n$ is unique, then so is $\pi^\ast_n$.
\end{thm}
{\it Proof: }
Let  
$Y^\ast_n$ denote the right hand side of \eqref{yast}. Then,
$\{(Y^\ast_n,Z^\ast_n)\}$ is the solution of the BS$\Delta$E
\begin{equation}\label{bsdeast}
    \Delta Y^\ast_n = - g_n(Z_n^\dagger) + (Z^\ast_n)^\top \Delta X_n,  \ \ Y^\ast_N = H
\end{equation}
by Theorem~\ref{thm:mart}.
 Note that
\begin{equation*}
    \Delta Y^\pi_n = -g_n(Z_n) + Z_n^\top \Delta X_n - \pi_n^\top \Delta X_n,\ \ Y^\pi_N = H
\end{equation*}
for a predictable process $\{Z_n\}$.
This means $\{(Y^\pi_n,Z^\pi_n)\}$, $Z^\pi_n = Z_n -\pi_n$ solves
the BS$\Delta$E
\begin{equation*}\label{bsde2}
\Delta Y^\pi_n = -h^\pi_n(Z^\pi_n) + (Z_n^\pi)^\top \Delta X_n,\ \ Y^\pi_N = H,
\end{equation*}
where $h^\pi_n(z) = g_n(z+ \pi_n)$.
For any predictable process $\pi$, we have 
$h^\pi_n \leq g_n(Z^\dagger_n)$.
Therefore $Y^\pi_n \leq Y^\ast_n$ by Theorem~\ref{thm:comp}.
Notice also that 
by choosing $\pi^\ast_n = Z^\dagger_n -Z^\ast_n$ we have
$g_n(Z^\dagger) = h_n^{\pi^\ast}(Z^\ast_n)$ so that
 $\{(Y^\ast_n,Z^\ast_n)\}$ satisfies the same BS$\Delta$E as $\{(Y^{\pi^\ast}_n,Z^{\pi^\ast}_n)\}$.
Hence, $Y^\ast_n = Y^{\pi^\ast}_n$.
\hfill$\square$\\

\begin{rem}\label{rem:hedge}
    \upshape
    Applying Theorem~\ref{thm:mart} to $g_n = 0$ and $Y_N = H$, we have
    \begin{equation*}
        H = \mathsf{E}_\mathsf{Q}[H] + \sum_{n=1}^N Z^H_n \Delta X_n
    \end{equation*}
    with $\{Z^H_n\}$ defined by \eqref{rep}.
    Therefore, the optimal strategy $\pi^\ast = Z^\dagger_n - Z^\ast_n = -Z^H_n + Z^\dagger_n - Z^g_n$ of Theorem~\ref{thm:max} is decomposed into the hedging part $-Z^H_n$ and the optimal investment part
    $Z^\dagger_n - Z^g_n$.
\end{rem}

\begin{ex}[ Locally entropic monetary utility]
    \upshape
    Consider \eqref{ex3g}. 
       Let $\Delta_d^\circ$ denote the interior of $\Delta_d$ and
      assume $\{(\hat{P}_{n,1},\dots, \hat{P}_{n,d})^\top\}$ to be $\Delta_d^\circ$-valued.
    Then, the map $z \mapsto g_n(z)$ is strictly concave and
   its unique maximizer is given by
   $Z^\dagger_n = z(y)$ of \eqref{syseq} for
   \begin{equation*}
       y = \frac{1}{G_n} \log \hat{P}_n
       :=  \frac{1}{G_n} (
       \log \hat{P}_{n,1}, \dots,
         \log \hat{P}_{n,d})^\top,
   \end{equation*}
   or equivalently (see \eqref{syseqsol}),
   \begin{equation}\label{ex5zd}
       Z^\dagger_n = \frac{1}{G_n} (\mathbf{v}\mathbf{v}^\top)^{-1}\mathbf{v}
       \log \hat{P}_n.
   \end{equation}
   Indeed, $v_j^\top Z^\dagger_n = G_n^{-1}\log \hat{P}_{n,j} - a(y)$ implies
   \begin{equation*}
       \sum_{j=0}^d v_j e^{-G_n v_j^\top Z^\dagger_n} \hat{P}_{n,j}= 0 
   \end{equation*}
   and so $\nabla g_n(Z^\dagger_n) =0$ for all $n$. We also have
   \begin{equation*}
   \begin{split}
      \frac{1}{G_n}\log B_n +    g_n(Z^\dagger_n) &= -\frac{1}{G_n} \log  \left(
         \sum_{j=0}^d e^{-G_n v_j^\top  Z^\dagger_n} \hat{P}_{n,j}\right) = 
          -\frac{1}{G_n} \log  \left(
          (d+1) e^{\frac{1}{d+1}\mathbf{1}^\top \log \hat{P}_n }
          \right)
          \\ &= \frac{1}{G_n} \frac{1}{d+1}\sum_{j=0}^d \log \frac{1}{(d+1)\hat{P}_{n,j}} = \frac{1}{G_n}D_{\mathrm{KL}}(Q_n ||\hat{P}_n) \geq 0,
   \end{split}
   \end{equation*}
   where $D_{\mathrm{KL}}$ denotes the Kullback-Leibler divergence on
   $\Delta_d$ and $Q_n = \mathbf{1}/(d+1)$.
   By \eqref{rep} and \eqref{ex5zd},
   \begin{equation}\label{exequi}
       \pi^\ast_n = Z^\dagger_n - Z^g_n - Z^H_n= - Z^H_n +  (\mathbf{v}\mathbf{v}^\top)^{-1}\mathbf{v}
       \left( \frac{\log \hat{P}_n}{G_n} - \hat{Y}_n \right), 
   \end{equation}
   where $\hat{Y}_n = (\hat{Y}_{n,0},\dots, \hat{Y}_{n,d})^\top$ and
   \begin{equation*}
       \hat{Y}_{n,j} =  \sum_{i=n+1}^n  \mathsf{E}_\mathsf{Q}[ g_i(Z^\dagger_i) |\mathscr{F}_{n-1}, \Delta X_n = v_j].
   \end{equation*}
   The first term $-Z^H_n$ is the hedging term as noted in Remark~\ref{rem:hedge}.
   The term
   \begin{equation*}
    Z^\dagger_n =    (\mathbf{v}\mathbf{v}^\top)^{-1}\mathbf{v}
       \frac{\log \hat{P}_n}{G_n} \approx 
        (\mathbf{v}\mathbf{v}^\top)^{-1}\mathbf{v}
       \frac{\hat{P}_n - \mathbf{1}}{G_n}
       = \frac{1}{G_n}  (\mathbf{v}\mathbf{v}^\top)^{-1}\hat{\mathsf{E}}[\Delta X_n | \mathscr{F}_{n-1}]
   \end{equation*}
   can be interpreted as the discrete counterpart of the Merton portfolio  (see e.g., Remark 8.9 of \cite{KS}). 
The term $\hat{Y}$ adjusts the expected return depending on the stochastic dynamics of
$D_{\mathrm{KL}}(Q_i ||\hat{P}_i)$ for $i \geq n+1$.
When $B_n = 1$ and $G_n = \gamma$  for all $n$ for a constant $\gamma > 0$ as in \eqref{ex1g}, we have
\begin{equation*}
\begin{split}
    \sum_{i=n+1}^n \mathsf{E}_\mathsf{Q}[g_i(Z^\dagger_i) | \mathscr{F}_n]
    &= \frac{1}{\gamma}
     \sum_{i=n+1}^n \mathsf{E}_\mathsf{Q}[D_\mathrm{KL}(Q_i||\hat{P}_i) | \mathscr{F}_n]
     \\ &= \frac{1}{\gamma}
      \mathsf{E}_\mathsf{Q}\left[ \log \frac{\mathrm{d}\mathsf{Q}}{\mathrm{d}\hat{\mathsf{P}}}
      - \log \mathsf{E}_\mathsf{Q}\left[  \frac{\mathrm{d}\mathsf{Q}}{\mathrm{d}\hat{\mathsf{P}}}\Bigg| \mathscr{F}_n\right]
      \Bigg| \mathscr{F}_n\right],
\end{split}    
\end{equation*}
where $\hat{\mathsf{P}}$ is associated with $\{\hat{P}_n\}$ by \eqref{phat}.
Note also that $Z^g_n = 0$ if $\{\hat{P}_n\}$, $\{G_n\}$ and $\{B_n\}$ are deterministic.
\end{ex}

\begin{ex}\label{VarSwap}
    \upshape
    Let $d=2$ and
    \begin{equation*}
        \mathbf{v} = \begin{pmatrix}
            0 & 1 & -1 \\
            -2c & c & c
        \end{pmatrix},
    \end{equation*}
    where $c \in \mathbb{R}$.
Then, we have
\begin{equation*}
    X_{n,2} = c\left(3 \sum_{k=1}^n |\Delta X_{k,1}|^2 - 2n\right).
\end{equation*}
Indeed, if $X_n = (n-a-b)v_0 + av_1 + bv_2$ for $(a,b) \in \mathbb{N}^2$, then
\begin{equation*}
    a + b = \sum_{k=1}^n |\Delta X_{k,1}|^2, \ \ 
    X_{n,2} = -2c(n-a-b) + c(a+b) = c(3(a+b)-2n).
\end{equation*}
Regarding $\{X_{n,1}\}$ as a price process of an asset,
the above identity allows us to interpret $X_{N,2}$ as an affine transform of
the variance swap payoff of the asset.
Further regarding $\mathsf{Q}$ as the pricing measure,
$\{X_{n,2}\}$ corresponds to the price process of the variance swap payoff.
Note that $\{X_{n,1}\}$ describes a trinomial model for the asset, which is not complete. 
The variance swap trading makes the $2$-dimensional market $\{X_n\}$ complete. 
\end{ex}
   \section{General equilibrium}
   \subsection{Formulation}
   We consider $m$ agents who 
   maximize respective utilities $\mathcal{E}_0^{(i)}(H^{(i)} + W_N(0,\pi^{(i)}))$, $i=1,\dots, m$ through trading strategies $\pi^{(i)}$ of the $d$-dimensional 
   asset $\{X_n\}$, where  $\mathcal{E}^{(i)}$ is the solution map of the BS$\Delta$E \eqref{bsde} with $g = g^{(i)} \in \mathcal{C}$,
   and $H^{(i)}$ is an  $\mathscr{F}_N$-measurable random variable representing
   an endowment for the agent $i$, 
    $i=1,\dots,m$.
Let $H_n$ denote the total supply vector of the asset vector $X_n$ and assume $\{H_n\}$ to be an $\mathbb{R}^d$-valued predictable process.
We say the market is in general equilibrium if there exist predictable processes $\pi^{(i)} = \{\pi^{(i)}_n\}$,
$i=1,\dots, m$ such that
\begin{enumerate}
    \item $\mathcal{E}_0^{(i)}(H^{(i)} + W_N(0,\pi^{(i)}))
    = \max_\pi \mathcal{E}_0^{(i)}(H^{(i)} + W_N(0,\pi))$ for all $i=1,\dots, m$, and
    \item $\sum_{i=1}^m \pi^{(i)}_n = H_n $ for all $n=1,\dots, N$,
\end{enumerate}
where the maximum is among all predictable processes $\pi$.

\begin{thm}\label{thm:supply}
    Let $H^{(i)}$ and $\tilde{H}^{(i)}$, $i=1,\dots,m$ be $\mathscr{F}_N$-measurable random variables, and $\{H_n\}$ and $\{\tilde{H}_n\}$ be $\mathbb{R}^d$-valued predictable processes
    satisfying
    \begin{equation*}
        \sum_{i=1}^m H^{(i)} + \sum_{n=1}^N H_n^\top \Delta X_n
        =
        \sum_{i=1}^m \tilde{H}^{(i)} + \sum_{n=1}^N \tilde{H}_n^\top \Delta X_n.
    \end{equation*}
    Then, 
the market with the endowments $H^{(i)}$ and the total supply $\{H_n\}$ is in general equilibrium if and only if so is
the one with the endowments $\tilde{H}^{(i)}$ and the total supply $\{\tilde{H}_n\}$.
\end{thm}
{\it Proof: }
Let $Z^{\dagger(i)}$, $Z^{H(i)}$ and  $Z^{g(i)}$ denote respectively
 $Z^{\dagger}$, $Z^{H}$ and  $Z^{g}$ in Theorem~\ref{thm:max}
with $g = g^{(i)}$ and $H = H^{(i)}$.
Then, by Theorem~\ref{thm:max},
\begin{equation*}
-H_n + \sum_{i=1}^m \pi^{(i)}_n = 
    -H_n - \sum_{i=1}^m Z^{H(i)} + \sum_{i=1}^m (Z^{\dagger(i)} - Z^{g(i)})
    =- Z^H + \sum_{i=1}^m (Z^{\dagger(i)} - Z^{g(i)}),
\end{equation*}
where $Z^H$ is defined by \eqref{rep} with
\begin{equation*}
    H =   \sum_{i=1}^m H^{(i)} + \sum_{n=1}^N H_n^\top \Delta X_n.
\end{equation*}
Therefore, whether $\sum_{i=1}^m \pi^{(i)}_n =H_n$ or not depends on $H^{(i)}$ and $\{H_n\}$ only  through $H$.
\hfill{$\square$}\\

 We are interested in conditions on $g^{(i)}$ for the market to be in general equilibrium, and so, assume hereafter $H^{(i)}=0$ ($i \geq 2)$ and $H_n =0$ ($n\geq 1$), without loss of generality in light of Theorem~\ref{thm:supply}.
 Let $H$ denote $H^{(1)}$.
 
   \subsection{Equilibrium for a single agent}
   Here we consider the case $m=1$.
   Put $g_n = g^{(1)}_n$ and $\pi^\ast_n = \pi^{(1)}_n$.
   The market is in general equilibrium if and only if
    $\mathcal{E}_0^g (H)
    = \max_\pi \mathcal{E}_0^{g}(H + W_N(0,\pi))$.

\begin{thm}\label{thm:sing1}
The market is in general equilibrium if \eqref{dagger} holds with
\begin{equation}\label{dagger2}
    Z^\dagger_n = (d+1)(\mathbf{v}\mathbf{v}^\top)^{-1}
       \mathsf{E}_\mathsf{Q}[\mathcal{E}^g_n(H)\Delta X_n | \mathscr{F}_{n-1}], \ \ n=1,\dots, N.
\end{equation}
\end{thm}
{\it Proof: }
By Theorem~\ref{thm:mart}, $(Y_n,Z_n) := (\mathcal{E}^g_n(H), Z_n^\dagger)$
solves \eqref{bsde} with $Y_N = H$.
Let $\{(Y^\ast_n,Z^\ast_n)\}$ be the solution of \eqref{bsdeast}.
Under \eqref{dagger}, 
since $Z_n = Z^\dagger_n$, 
$\{(Y_n,Z_n)\}$ solves \eqref{bsdeast} as well, which means
$(Y_n,Z_n) = (Y_n^\ast,Z_n^\ast)$.
Thus we have $Z_n^\dagger = Z^\ast_n$.
 By Theorem~\ref{thm:max}, we conclude $\pi^\ast_n = 0$ and so, the market is in general equilibrium.
\hfill{$\square$}

 \begin{thm}\label{thm:sing2}
Assume that there uniquely exists $\{Z^\dagger_n\}$ for \eqref{dagger} to hold.
If the market is in general equilibrium, then we have \eqref{dagger2}.
\end{thm}
{\it Proof: }
Let $\{(Y^\ast_n,Z^\ast_n)\}$ be the solution of \eqref{bsdeast}.
       By Theorem~\ref{thm:max},  for the market to be in general equilibrium we need $Z^\dagger_n = Z^\ast_n$.
Because of this identity, $\{(Y^\ast_n,Z^\ast_n)\}$ solves
\eqref{bsde} with $Y_N = H$, meaning that 
$Y^\ast_n = \mathcal{E}^g_n(H)$.
We then have \eqref{dagger2} by Theorem~\ref{thm:mart}.
\hfill{$\square$}

\begin{rem}\upshape \label{rem:rec}
    Although \eqref{dagger2} (with \eqref{dagger}) is a fixed point type equation for $g$,
    it is actually a recurrence relation for the sequence $\{g_n\}$. Indeed,
    $\mathcal{E}^g_n(H)$ is determined without specifying $g_n$ in the sense that 
    $\mathcal{E}^g_n(H) = \mathcal{E}^{\hat{g}}_n(H)$
if $g_k = \hat{g}_k$ for all $k > n$.
\end{rem}
\begin{thm}\label{thm:sing}
   Let 
   $f_n:\Omega \times \mathbb{R}^d \times \Delta_d \to \mathbb{R}$ be $\mathscr{F}_{n-1}\otimes \mathscr{B}(\mathbb{R}^d \times \Delta_d)$-measurable, concave on $\mathbb{R}^d$ and continuously differentiable on $\mathbb{R}^d$ with
   $\nabla f_n$ taking values in $\Theta$, and 
   $0 \in \nabla f_n(z,\Delta_d)$ for all $z \in \mathbb{R}^d$,
   $n=1,\dots, N$.
   Then,
   there exists a  $\Delta_d$-valued predictable process $\{\hat{P}_n\}$ such that the market with $g =\{g_n\}$,
    $g_n(z) = f_n(z,\hat{P}_n)$,  is in general equilibrium.
\end{thm}
{\it Proof: }
 Note first that $g \in \mathcal{C}$ by Remark~\ref{rem:comp}.
We construct $\{\hat{P}_n\}$ inductively.
By the assumption, there exists $\hat{P}_N$ such that
$\nabla f_N(Z^\dagger_N,\hat{P}_N) = 0$ 
for $Z^\dagger_N$ defined by \eqref{dagger2} for $n=N$ and $\mathcal{E}^g_N(H) = H$.
  Given $\hat{P}_k$ for $k \geq n+1 $,
  let $Z^\dagger_n$ be defined by \eqref{dagger2} with
  $g_k(z) = f_k(z,\hat{P}_k)$ (see Remark~\ref{rem:rec}).
  By the assumption, there exists $\hat{P}_n$ such that
$\nabla f_n(Z^\dagger_n,\hat{P}_n) = 0$.
Since $z \mapsto f_n(z,\hat{P}_n)$ is concave,
  $Z^\dagger_n$ is a maximizer of $g_n(z) = f_n(z,\hat{P}_n)$.
  The result then follows from Theorem~\ref{thm:sing1}.
\hfill{$\square$}

\begin{ex}[Locally entropic monetary utility] \label{Ex41}
    \upshape
    We consider \eqref{ex3g} again  with
     $\{\hat{P}_n\}$,
     $\hat{P}_n = (\hat{P}_{n,1},\dots, \hat{P}_{n,d})^\top$, being $\Delta_d^\circ$-valued.
    The driver $g_n$ is of the form 
    $g_n(z) = (f(G_nz,\hat{P}_n) -\log B_n)/G_n $, where
    \begin{equation}\label{fzp}
        f(z,p) = -\log \sum_{j=0}^d e^{-z^\top v_j}p_j,\ \ z \in \mathbb{R}^d, \ \ p = (p_0,\dots,p_d)^\top \in \Delta_d.
    \end{equation}
    By \eqref{ex5zd},
    we have $\nabla f_n(Z_n,\hat{P}_n) = 0$ if and only if
    \begin{equation*}
        Z_n  = \frac{1}{G_n} (\mathbf{v}\mathbf{v}^\top)^{-1}\mathbf{v}
       \log \hat{P}_n.
    \end{equation*}
Since $\mathbf{v}$ has rank $d$ and $\mathbf{v}\mathbf{1} = 0$,
by Theorems~\ref{thm:sing1} and \ref{thm:sing2},
the market is in general equilibrium if and only if
   \begin{equation}\label{ex5eq}
       \hat{P}_{n,j} = \frac{e^{G_n\bar{Y}_{n,j}}}{\sum_{k=0}^de^{G_n\bar{Y}_{n,k}} }, \ \ 
       \bar{Y}_{n,j} = \mathsf{E}_\mathsf{Q}[\mathcal{E}^g_n(H) | \mathscr{F}_{n-1},\Delta X_n = v_j]
       \ \ j=0,\dots, d.
   \end{equation}
   Note that
   \begin{equation*}
       \hat{P}_{n,j}1_{\{\Delta X_n = v_j\}} = \frac{e^{G_n\mathcal{E}^g_n(H)}}{(d+1)\mathsf{E}_\mathsf{Q}[e^{G_n\mathcal{E}^g_n(H)}|\mathscr{F}_{n-1}]}1_{\{\Delta X_n=v_j\}},\ \ j=0,\dots,d
   \end{equation*}
   under \eqref{ex5eq}.
   Therefore, 
   the market is in general equilibrium
   if and only if
   \begin{equation*}
       \frac{\mathrm{d}\hat{\mathsf{P}}}{\mathrm{d}\mathsf{Q}}=\prod_{n=1}^N\frac{e^{G_n\mathcal{E}^g_n(H)}}{\mathsf{E}_{\mathsf{Q}}[e^{G_n\mathcal{E}^g_n(H)}|\mathscr{F}_{n-1}]},
   \end{equation*}
   where $\hat{\mathsf{P}}$ is defined by \eqref{phat}.
 Further by \eqref{trans}, we have
 \begin{equation*}
 -\frac{1}{G_n}\log \hat{E}[
      e^{-G_n \mathcal{E}^g_n(H)} |\mathscr{F}_{n-1}
      ]  = 
     \mathcal{E}^g_{n-1}(H) + \frac{1}{G_n} \log B_n.
 \end{equation*}
   This implies
   \begin{equation}\label{ex41}
       \frac{\mathrm{d}\mathsf{Q}}{\mathrm{d}\hat{\mathsf{P}}}=\prod_{n=1}^N\frac{e^{-G_nY^\ast_n}}{\hat{\mathsf{E}}[e^{-G_n Y^\ast_n}|\mathscr{F}_{n-1}]}=
       \exp\left\{
\sum_{n=1}^N -G_n (\mathcal{E}^g_n(H)-\mathcal{E}^g_{n-1}(H))
       \right\} \prod_{n=1}^N B_n. 
   \end{equation}
   When $B_n = 1$ and $G_n = \gamma$  for all $n$ for a constant $\gamma > 0$ as in \eqref{ex1g}, we have
   \begin{equation}\label{ess}
     \frac{\mathrm{d}\mathsf{Q}}{\mathrm{d}\hat{\mathsf{P}}}= 
     \frac{e^{-\gamma H}}{\hat{\mathsf{E}}[e^{-\gamma H}]},
\end{equation}
which characterizes the equilibrium probability measure $\hat{\mathsf{P}}$.
\end{ex}

\begin{rem}
    \upshape
    In addition to the conditions of Theorem~\ref{thm:sing},
    if $f_n$ is of the form
    $f_n(z,p) = f(G_n z,p)/G_n$ 
    for a smooth function $f:\mathbb{R}^d \times \Delta_d \to \mathbb{R}$
    with $\nabla f(0,p) = p$
    and a positive predictable process $\{G_n\}$
    as in Example~\ref{Ex41} with $B_n = 1$, then we have
    \begin{equation*}
        g_n(z) = f_n(z,\hat{P}_n) \approx z^\top \hat{P}_n
    \end{equation*}
    considering $G_n$ to be small.
    This approximation implies in turn
    \begin{equation*}
        \mathcal{E}^g_n(Y) \approx \hat{\mathsf{E}}[Y |\mathscr{F}_n]
    \end{equation*}
    in light of Theorem~\ref{thm:linear},
     where $\hat{\mathsf{E}}$ is the expectation under $\hat{\mathsf{P}}$  defined by \eqref{phat}.
     Therefore in this case, extending
 Example~\ref{Ex41}, we can interpret  $G_n$ as a risk aversion parameter
    and $\hat{\mathsf{P}}$ as the belief of the agent.
    If the market is in general equilibrium,
   \begin{equation*}
       \hat{\mathsf{E}}[\Delta X_n | \mathscr{F}_{n-1}] = \mathbf{v}\hat{P}_n
   \end{equation*}
   is interpreted as an equilibrium return.
\end{rem}

\begin{ex}\label{ex:markov}
    \upshape
    If $H = h(X_N)$ and $g_n(z) = f_n(X_{n-1},z)$ for deterministic functions
    $h$ and $f_n$ as in Section~\ref{FK} and if $f(x,z)$ is strictly concave and continuously differentiable in $z$, then the condition \eqref{dagger2} 
    follows from the deterministic identity
    \begin{equation*}
        \nabla_z f_n(x,(\mathbf{vv}^\top)^{-1}\mathbf{v}\mathcal{N}u_n(x)) = 0,\ \ x \in L, \ \  n=1,\dots, N,
    \end{equation*}
    where $\mathcal{N}u_n$ is as in Section~\ref{FK}.
    For example, for the market with no random endowments ($H^{(i)}=0$) and 
    unit total supply ($H_n = \mathbf{1}$), we have $H = \mathbf{1}^\top X_N$.
\end{ex}

\begin{ex}
    \upshape
    Consider \eqref{ex3g} with $B_n = 1$, $G_n = \gamma_n(X_{n-1})$ and 
    $\hat{P}_n = p_n(X_{n-1})$ for deterministic functions $\gamma_n:L \to (0,\infty)$ and $p_n: L \to \Delta_d^\circ$.
    Assume also $H  = h (X_N)$ as in Example~\ref{ex:markov}.
    Then, from Examples~\ref{Ex41} and \ref{ex:markov}, the market is in general equilibrium if
    \begin{equation*}
        \mathbf{v}\log p_n(x) = \gamma_n(x) \mathbf{v} \mathcal{N}u_n(x), \ \ x \in L,
        \ \ n=1,\dots, N.
    \end{equation*}
    The function $u_n(x)$ is computed backward inductively without using $p_n(x)$.
    For a given function $\gamma_n(x)$, there exists a unique $p_n(x) \in \Delta_d^\circ$ satisfying this equation for each $x \in L$. 
    For the sequence of such functions $p_n(x)$ obtained in the backward manner, the $\Delta_d$-valued sequence $\hat{P}_n = p_n(X_{n-1})$ defines a unique equilibrium probability $\hat{\mathsf{P}}$ by \eqref{phat} associated with the sequence $G_n = \gamma_n(X_{n-1})$.
    The equilibrium return is approximated as 
    \begin{equation*}
        \hat{\mathsf{E}}[\Delta X_n | \mathscr{F}_{n-1}] =
        \mathbf{v} \hat{P}_n\approx 
        \mathbf{v}\log \hat{P}_n =  
        \gamma_n(X_{n-1})\mathbf{v} \mathcal{N}u_n(X_{n-1})
    \end{equation*}
    using $\mathbf{v1} = 0$.
\end{ex}

\subsection{The representative agent}
Here we consider the case $m>1$.
For functions  $f^{(1)}$ and $f^{(2)}$ on $\mathbb{R}^d$, define
the sup-convolution $f^{(1)} \square f^{(2)}$ by
\begin{equation*}
    f^{(1)}\square f^{(2)}(z) = \sup_{x \in \mathbb{R}^d} \{f^{(1)}(x) + f^{(2)}(z-x)\}.
\end{equation*}
For the drivers $g^{(i)} = \{g^{(i)}_n\}$, $i=1,\dots,m$ of the $m$ agents' utilities, let
\begin{equation}\label{representative}
    g_n(z) = g^{(1)}_n \square \dots \square g^{(m)}_n(z).
\end{equation}
A single agent whose utility is $\mathcal{E}^g_0$ with $g = \{g_n\}$ defined by \eqref{representative} and whose endowment is $H$ is called the representative agent of the market.

\begin{lem} 
For all $n$,
    \begin{equation}\label{supsum}
    \sup_{z \in \mathbb{R}^d} g_n(z)
    = \sum_{i=1}^m  \sup_{z \in \mathbb{R}^d} g_n^{(i)}(z).
\end{equation}
\end{lem}
{\it Proof: }
It is trivial that $g_n(z)$ is upper bounded by the right hand side sum for any $z$ hence its supremum. Conversely, let $\{z^{(i)}_k\}$ be a such sequence that of $ \lim_{k \to \infty} g_n^{(i)}(z^{(i)}_k) =  \sup_{z \in \mathbb{R}^d} g_n^{(i)}(z)$ for each $i$. Then,
\begin{equation*}
    \sum_{i=1}^m g_n^{(i)}(z^{(i)}_k) \leq 
    g_n\left(\sum_{i=1}^m z^{(i)}_k \right) \leq  \sup_{z \in \mathbb{R}^d} g_n(z)
\end{equation*}
for any $k$, hence its limit.
\hfill{$\square$}

\begin{thm}\label{thm:rep}
Assume that there exist
unique predictable processes $\{Z^{\dagger(i)}_n\}$ for each $i=1,\dots,m$ such that
\begin{equation*}
   g^{(i)}_n(Z^{\dagger(i)}_n) = \sup_{z \in \mathbb{R}^d}g^{(i)}_n(z)
\end{equation*}
for all $n = 1,\dots, N$. 
Assume further that the maximizer of $z \mapsto g_n(z)(\omega)$ is also unique for all $n$ and all $\omega \in \Omega$.
  The market for the $m$ agents is in general equilibrium 
  if and only if so is the one for the representative agent.
\end{thm}
{\it Proof: }
Let  $\{(Y^{\ast(i)}_n,Z^{\ast(i)}_n)\}$ be the solution of
\begin{equation*}
    \Delta Y^{\ast (i)}_n = -g^{(i)}_n(Z^{\dagger (i)}_n) + (Z^{\ast (i)}_n)^\top \Delta X_n
\end{equation*}
with $Y^{\ast (1)}_N = H$ and  $Y^{\ast (i)}_N = 0$ for $i \geq 2$.
By Theorem~\ref{thm:max}, the unique optimal strategy $\pi^{(i)}$ for the agent $i$ is given by $\pi^{(i)}_n = Z^{\dagger (i)}_n -Z^{\ast(i)}_n $.
 By \eqref{supsum}, we have
\begin{equation*}
   \sum_{i=1}^m  g^{(i)}_n(Z^{\dagger(i)}_n)
  = \sup_{z \in \mathbb{R}^d}g_n(z) =   g_n(Z^\dagger_n),
    \ \ 
     Z^\dagger_n := \sum_{i=1}^m Z^{\dagger(i)}_n.
\end{equation*}
Therefore, 
\begin{equation*}
    Y^\ast_n := \sum_{i=1}^m Y^{\ast(i)}_n, \ \ 
    Z^\ast_n :=  \sum_{i=1}^m Z^{\ast (i)}_n
\end{equation*}
solves
\begin{equation*}
    \Delta Y^\ast_n = -g_n(Z^\dagger_n) + (Z^\ast_n)^\top \Delta X_n, \ \ Y^\ast_N = H.
\end{equation*}
Since $Z^\dagger_n$ is the unique maximizer of $g_n$,
the unique optimal strategy for the representative agent is
$\pi^\ast_n = Z^\dagger_n - Z^\ast_n$ again by Theorem~\ref{thm:max}.
Hence,
\begin{equation*}
    \sum_{i=1}^m \pi^i_n =  \sum_{i=1}^m
   Z^{\ast (i)}_n - \sum_{i=1}^m Z^{\dagger (i)}_n = Z^\ast_n - Z^\dagger_n = \pi^\ast_n.
\end{equation*}
Therefore $\sum_{i=1}^m \pi^i_n = 0$ if and only if $\pi^\ast_n = 0$.
\hfill{$\square$}

\subsection{Equilibrium for multiple agents with a common belief}
Here we suppose
$g^{(i)}_n(z) = f_n(G^{(i)}_nz,\hat{P}_n)/G^{(i)}_n$, where
$f_n: \Omega \times \mathbb{R}^d \times \Delta_d \to \mathbb{R}$
and $\hat{P}_n$ are as in Theorem~\ref{thm:sing}, 
and $\{G^{(i)}_n\}$ is a positive predictable process for each $i=1,\dots, m$.
By induction, we can show that
\begin{equation*}
     g_n(z) := g^{(1)}_n \square \dots \square g^{(m)}_n(z)
     = \frac{1}{G_n}  f_n(G_nz,\hat{P}_n)
\end{equation*}
with
\begin{equation}\label{agr}
    \frac{1}{G_n} = \sum_{i=1}^m  \frac{1}{G^{(i)}_n}.
\end{equation}
\begin{thm}
There exists a $\Delta_d$-valued predictable process $\{\hat{P}_n\}$ such that
the market is in general equilibrium.
\end{thm}
{\it Proof: } It follows from Theorems~\ref{thm:sing} and \ref{thm:rep}.
\hfill{$\square$}
\begin{ex}[Locally entropic monetary utility]
    \upshape
    Let $f_n(z,p) = f(z,p)$ defined by \eqref{fzp}.
  The predictable process  $\{G^{(i)}_n\}$ and $\{G_n\}$ are then interpreted as the local risk aversion parameter for the agent $i$ and for the representative agent respectively.
    Then, \eqref{ex5eq} defines the unique sequence
    $\{\hat{P}_n\}$ such that
the market is in general equilibrium.
The equilibrium probability measure $\hat{\mathsf{P}}$ satisfies \eqref{ex41} with $B_n = 1$.
When $G_n = \gamma$ for all $n$ for a constant $\gamma > 0$,
then the representative agent has the exponential utility \eqref{ex1}
and the equilibrium probability measure $\hat{\mathsf{P}}$ is characterized by 
\eqref{ess}.
\end{ex}
\subsection{Equilibrium under heterogeneous beliefs}
Here we 
assume locally entropic monetary utilities
$g^{(i)}_n(z) = f(G^{(i)}_nz,\hat{P}^{(i)}_n)/G^{(i)}_n$,
where $f$ is defined by \eqref{fzp} and $\{G^{(i)}_n\}$ and
$$\{\hat{P}^{(i)}_n\} = 
\left\{\left(\hat{P}^{(i)}_{n,0},\dots,\hat{P}^{(i)}_{n,d}\right)^\top \right\} $$ are respectively
$(0,\infty)$-valued and $\Delta_d^\circ$-valued predictable processes for each $i=1,\dots,m$.
Each sequence $\{\hat{P}^{(i)}_n\}$ determines a probability measure $\hat{P}^{(i)}$ on $\mathscr{F}_N$ by \eqref{phat}, which is interpreted as the agent $i$'s belief on the law of $\{X_n\}$.
\begin{thm} \label{thm:het}
Define $\{G_n\}$ by \eqref{agr}. Then,
    \begin{equation*}
        g^{(1)}_n \square \dots \square g^{(m)}_n(z) = 
        \frac{1}{G_n}f(G_nz,\tilde{P}_n) - \frac{1}{G_n}\log C_n,
    \end{equation*}
    where $ \tilde{P}_n = ( \tilde{P}_{n,0},\dots,  \tilde{P}_{n,d})^\top$ and
    \begin{equation*}
        \tilde{P}_{n,j} = \frac{1}{C_n}\prod_{i=1}^m (\hat{P}^{(i)}_{n,j})^{G_n/G_n^{(i)}}, \ \ 
        C_n = \sum_{j=0}^d \prod_{i=1}^m (\hat{P}^{(i)}_{n,j})^{G_n/G_n^{(i)}}.
    \end{equation*}
\end{thm}
{\it Proof: }
The case $m=2$ follows by solving the equation
\begin{equation*}
    \nabla f(G^{(1)}_nx,\hat{P}^{(1)}_n) = \nabla f(G^{(2)}_n(z-x),\hat{P}^{(2)}_n)
\end{equation*}
in $x \in \mathbb{R}^d$; see Lemma~\ref{lemA1}. The general case then follows by induction.
\hfill{$\square$}

\begin{rem}
By Lemma~\ref{lemA2}, we have
    $C_n \leq 1$ with equality holding if and only if $\hat{P}^{(i)}_n =\hat{P}^{(1)}_n$ for all $i$. 
\end{rem}

By Theorem~\ref{thm:het},
the representative agent's market falls into Example~\ref{Ex41}.
In particular, when $G_n = \gamma>0$ (a constant),
the market is in general equilibrium if and only if
$\tilde{P}_n =\hat{P}_n$, where
 \begin{equation*}
       \frac{\mathrm{d}\mathsf{Q}}{\mathrm{d}\hat{\mathsf{P}}}=
       \exp\left\{
\sum_{n=1}^N -G_n (\mathcal{E}^g_n(H)-\mathcal{E}^g_{n-1}(H))
       \right\} \prod_{n=1}^N C_n = 
        \frac{e^{-\gamma H} \prod_{n=1}^N C_n }{\hat{\mathsf{E}}[e^{-\gamma H} \prod_{n=1}^N C_n ]}.
   \end{equation*}
To highlight the outcome of heterogeneous beliefs,
let us further assume there are only two agents ($m=2$) with constant risk aversion $G^{(i)}_n = \gamma_i>0$ and with no endowment ($H=0$).
In this case, $\nabla g^{(i)}(0,\mathsf{Q}) = 0$ and so,
if the two agents have a common belief, it must be 
the martingale measure $\mathsf{Q}$ for the market to be in general equilibrium. The optimal strategies are simply $\pi^{(1)}_n = \pi^{(2)}_n = 0$. On the other hand, for any $\Delta_d^\circ$-valued deterministic sequence
$\{\hat{P}^{(1)}_n\}$, by choosing $\hat{P}^{(2)}_n$ as
\begin{equation*}
\hat{P}^{(2)}_{n,j} = \frac{(\hat{P}^{(1)}_{n,j})^{-\gamma_2/\gamma_1}}{\sum_{k=0}^d (\hat{P}^{(1)}_{n,k})^{-\gamma_2/\gamma_1}},
\end{equation*}
we have $\tilde{P}_n = Q_n$ with $C_n$ being deterministic, which makes this market with heterogeneous beliefs be in general equilibrium.
The individual optimal strategies $\pi^{(1)}_n = - \pi^{(2)}_n$ are nonzero; the agents bet on their beliefs.

Another observation is that the equilibrium return is mostly affected by the belief of the least risk averse agent. Indeed, if 
$G^{(1)}_n \ll G^{(i)}_n$ for $i \geq 2$, we have
$G_n/G_n^{(1)} \approx 1$ while $G_n/G_n^{(i)} \approx 0$ for $i \geq 2$.
Therefore $\tilde{P}_n \approx \hat{P}^{(1)}_n$.

\begin{appendix}
\section{A computation of sup-convolution}
\begin{lem}\label{lemA1}
Let $\alpha > 0$, $\beta > 0$, 
$p_j > 0$, $q_j > 0$, $j=0,\dots,d$,
 and $z \in \mathbb{R}^d$. Then,
\begin{equation*}
\begin{split}
 &  \sup_{x \in \mathbb{R}^d}
\left\{
-\frac{1}{\alpha}\log \sum_{j=0}^d e^{-\alpha x^\top v_j}p_j
-\frac{1}{\beta}\log \sum_{j=0}^d e^{-\beta (z-x)^\top v_j}q_j
\right\} 
\\ &= -\frac{1}{\gamma} \log \sum_{j=0}^d e^{-\gamma z^\top v_j} p_j^{\gamma/\alpha}q_j^{\gamma/\beta},
\end{split}
\end{equation*}
where $\gamma = 1/(1/\alpha + 1/\beta)$.
\end{lem}
{\it Proof: }
The first order condition is
\begin{equation*}
    \sum_{j=0}^d v_j \left(\frac{e^{-\alpha x^\top v_j} p_j}{S(x,\alpha,\{p_k\})} - \frac{e^{-\beta (z-x)^\top v_j}q_j}{S(z-x,\beta,\{q_k\})}\right) = 0,
\end{equation*}
where
\begin{equation*}
    S(x,\alpha,\{p_k\}) = \sum_{j=0}^d e^{-\alpha x^\top v_j}p_j.
\end{equation*}
Since $\mathbf{v}\mathbf{1}= 0$ and $\mathrm{rank}\, \mathbf{v} = d$, 
the first order condition is met if and only if
\begin{equation*}
    -\alpha x^\top v_j + \log p_j = -\beta(z-x)^\top v_j + \log q_j + c(x),\ \ 
    j=0,\dots,d
\end{equation*}
for a function $c$. 
Substituting
\begin{equation*}
    x^\top v_j = \frac{\beta}{\alpha + \beta}z^\top v_j + \frac{\log p_j - \log q_j -c(x)}{\alpha + \beta},
\end{equation*}
we obtain 
\begin{equation*}
    \begin{split}
      &  -\frac{1}{\alpha}\log \sum_{j=0}^d e^{-\alpha x^\top v_j}p_j
      = -\frac{1}{\alpha}\log \sum_{j=0}^d e^{-\gamma z^\top v_j}p_j^{\gamma/\alpha}q_j^{\gamma/\beta} - \frac{c(x)}{\alpha + \beta}
\\
&-\frac{1}{\beta}\log \sum_{j=0}^d e^{-\beta (z-x)^\top v_j}q_j
  = -\frac{1}{\beta}\log \sum_{j=0}^d e^{-\gamma z^\top v_j}p_j^{\gamma/\alpha}q_j^{\gamma/\beta} + \frac{c(x)}{\alpha + \beta},
\end{split}
\end{equation*}
hence the result.\hfill{$\square$}

\begin{lem}\label{lemA2}
    Let $(p_{i,0},\dots,p_{i,d})^\top$,  $i=1,\dots,m$
    be $m$ points in $\Delta_d^\circ$.
    Let $\gamma_i > 0$ for $i=1,\dots,m$ and
    \begin{equation*}
        \gamma = \left(\sum_{i=1}^m \frac{1}{\gamma_i}\right)^{-1}.
    \end{equation*}
    Then,
    \begin{equation*}
        \sum_{j=0}^d \prod_{i=1}^m p_{i,j}^{\gamma/\gamma_i} \leq 1.
    \end{equation*}
\end{lem}
{\it Proof: } The case $m=1$ is trivial.
Let
\begin{equation*}
        \hat\gamma_k = \left(\sum_{i=1}^k \frac{1}{\gamma_i}\right)^{-1}, \ \ k=1,\dots,m.
    \end{equation*}
If the inequality is true when $m=k$, then
\begin{equation*}
\begin{split}
      \sum_{j=0}^d \prod_{i=1}^{k+1} p_{i,j}^{\hat\gamma_{k+1}/\gamma_i} 
      &= \sum_{j=0}^d p_{k+1,j} \left(\frac{\prod_{i=1}^{k} p_{i,j}^{\hat\gamma_k/\gamma_i} }{p_{k+1,j}}\right)^{\hat\gamma_{k+1}/\hat\gamma_k}
    \leq \left(  \sum_{j=0}^d \prod_{i=1}^{k} p_{i,j}^{\hat\gamma_{k}/\gamma_i} \right)^{\hat\gamma_{k+1}/\hat\gamma_k}
    \leq 1
\end{split}
\end{equation*}
by Jensen's inequality. We obtain the result by induction. \hfill{$\square$}
\end{appendix}

\end{document}